\documentclass{article}
\usepackage[margin=1in]{geometry}
\usepackage{amstext,amsfonts,amssymb,amscd,amsbsy,amsmath,verbatim,amsthm} 
\usepackage[all]{xy}
\usepackage{hyperref}
\usepackage{color,enumerate,float,graphicx}
\usepackage{chapterbib,appendix}
\usepackage[utf8]{inputenc} 

\renewcommand{\thesection}{\arabic{section}}

\numberwithin{equation}{section}
\newcounter{defi}
\renewcommand{\thedefi}{\arabic{section}.\arabic{defi}}
\numberwithin{defi}{section} 
\newcommand{\newdefi}[1]{\refstepcounter{defi}   \textbf{#1~\thedefi.}} 
\newcommand{\bewende}{\hfill $\square$\vspace{2mm}\\}

\newcommand{\tfo}{\tilde{F}_1}
\newcommand{\tft}{\tilde{F}_2}
\newcommand{\Hcurl}{H_{\mathrm{curl}}}

\newcommand{\wcal}{\mathcal{W}}
\newcommand{\vcal}{\mathcal{V}}

\newcommand{\pnu}{P(i\nu)}	
\newcommand{\pna}{P(i\nabla)}

\newcommand{\on}{~\mathrm{on}~}			
\newcommand{\inm}{~\mathrm{in}~}
\newcommand{\andm}{~\mathrm{and}~}
\newcommand{\diag}[1]{\mathrm{diag}\big(#1\big)}
\newcommand{\XX}{\mathbb{X}}
\newcommand{\RR}{\mathbb{R}}
\newcommand{\CC}{\mathbb{C}}

\newcommand{\td}{\mathrm{diag}\left(\tilde{\mu},\tilde{\gamma}\right)}		
\newcommand{\tdg}{\mathrm{diag}\left(\tilde{\gamma},\tilde{\mu}\right)}
\newcommand{\tdo}{\mathrm{diag}\left(\tilde{\mu},\overline{\tilde{\gamma}}\right)}			
\newcommand{\tdgo}{\mathrm{diag}\left(\overline{\tilde{\gamma}},\tilde{\mu}\right)}
\newcommand{\ifour}{\dot{I}_4}
\newcommand{\exi}{e^{i\xi\cdot x}}

\makeatletter
\def\blfootnote{\gdef\@thefnmark{}\@footnotetext}
\makeatother

\begin{document}

\title{Inverse problems for Maxwell's equations in a slab \\with partial boundary data\blfootnote{Date: \today}}
\author{Monika Pichler}
\date{}
\maketitle

\begin{abstract}
We consider two inverse boundary value problems for the time-harmonic Maxwell equations in an infinite slab. Assuming that tangential boundary data for the electric and magnetic fields at a fixed frequency is available either on subsets of one boundary hyperplane, or on subsets of different boundary hyperplanes, we show that the electromagnetic material parameters, the conductivity, electric permittivity, and magnetic permeability, are uniquely determined by these partial measurements.
\end{abstract}

\section{Introduction}

In this work, we investigate the unique determination of the electromagnetic material properties of an object from surface measurements. Inverse boundary value problems of this kind arise in many physical situations where one wishes to determine certain properties of a body from measurements taken on its surface, or at a distance from it. In many applications, the situation is described by a partial differential equation, and the inverse problem mathematically amounts to reconstructing parameters of the equation from boundary data of solutions to the equation. Since the seminal paper of A.P. Calder\'{o}n \cite{C1980} formulating the inverse problem for the conductivity equation on a bounded domain, and the subsequent work by Sylvester and Uhlmann \cite{SU1987} showing unique solvability of this problem for a smooth conductivity, many advances have been made in this field. The method developed in \cite{SU1987} of constructing special exponentially growing solutions to the equation, called \emph{complex geometrical optics} solutions (CGO solutions), has proved to be applicable to many different situations, and it remains the standard method for showing unique solvability of such inverse boundary value problems. A survey of advances made in the field can be found in \cite{U2009}.

This method is also the basis for the study of the inverse boundary value problem for Maxwell's equations. This problem was first formulated on a bounded domain in \cite{SIC1992}, where the linearization of the problem was investigated. Ola, P\"{a}iv\"{a}rinta, and Somersalo \cite{OPS1993} first proved global uniqueness of the solution to the inverse problem with full data for smooth parameters; an alternative proof was later given in \cite{OS1996}, by relating Maxwell's equations to a vector Schr\"{o}dinger equation. Caro and Zhou \cite{CZ2014} reduced the required regularity to continuously differentiable parameters, and the author recently showed that uniqueness also holds for Lipschitz parameters \cite{P2018}. 

In recent years, partial data inverse problems, in which the boundary data is available only on a subset of the boundary, have become a focus of attention. This is due to their practical importance, since in applications, part of the surface of the object of interest may not be available for measurements, or it may simply be too costly to perform measurements on the whole surface. Mathematically, the consequence of partial boundary data is that the integral identity that is the starting point for the uniqueness proof will involve boundary integrals with unknown functions. Two main methods have been found to be effective in dealing with such problems; on the one hand, the use of Carleman estimates to control the size of solutions on inaccessible parts of the boundary, on the other hand, reflection methods used to construct solutions that vanish on the inaccessible part, assuming this part of the boundary has a suitable geometry.

Carleman estimates were first used in a partial data inverse problem for the conductivity and Schr\"{o}dinger equations by Bukhgeim and Uhlmann \cite{BU2002}. 
Kenig, Sj\"{o}strand, and Uhlmann \cite{KSU2007} introduced Carleman estimates with nonlinear weight functions which allowed them to significantly improve on the previous result.
The method has since been extended to different inverse problems, such as for Schr\"{o}dinger operators with magnetic potential \cite{DKSU2007}, 
and for Maxwell's equations on a manifold \cite{COST2017}. 

Isakov \cite{I2007} introduced the reflection method and showed unique determination of the conductivity from Dirichlet and Neumann data given on some subset of the boundary, assuming that the inaccessible part of the boundary is part of a plane or of a sphere. This restriction allows to reflect solutions in such a way across the inaccessible part that one obtains solutions with vanishing boundary values. Caro, Ola, and Somersalo \cite{COS2009} applied the reflection argument in studying a partial data inverse problem for Maxwell's equations on a bounded domain. 

The geometry of a slab, an infinite domain bounded by two parallel hyperplanes, has been of considerable interest in view of applications in modeling waveguides and in medical imaging, among others. Consequently, many different types of inverse problems 
have been studied in this setting. The question of identifying an object embedded in a homogeneous slab was considered in \cite{I2001,SW2006}, inverse scattering problems in acoustics and optics have been studied for example in \cite{DM2006,ERY2008}, and \cite{MS2001} considered an inverse problem in optical diffusion tomography. 

There are also a number of recent results on partial data problems in a slab for the conductivity and 
Schr\"{o}dinger equations. Due to the geometry of the boundary, reflection arguments are a powerful tool in this scenario; however, Carleman estimates have been employed as well.
Li and Uhlmann \cite{LU2010} studied inverse problems for the conductivity and scalar Schr\"{o}dinger equations in two different cases, namely with partial Dirichlet and Neumann data given (i) on the same boundary hyperplane or (ii) on opposite boundary hyperplanes. 
The authors use the Carleman estimate derived in \cite{BU2002} 
to show that the boundary integrals are negligible in case (i); in case (ii), a reflection argument is employed. In \cite{KLU2012}, the inverse problem for a Schr\"{o}dinger operator with a magnetic potential was studied in the same setting,  
using reflection arguments to construct the necessary CGO solutions; in \cite{Li2012,Li2012opp} a matrix Schr\"{o}dinger operator was considered in each of the cases, using a reflection argument and a combination of Carleman estimates and reflection arguments, respectively. 

In this paper, we want to investigate each of these two scenarios in the partial data inverse problem for Maxwell's equations in the slab.  
We show unique solvability in these cases  
by relating them to the bounded domain setting investigated in \cite{COS2009}, and employing arguments similar to those used there, as well as arguments used in \cite{KLU2012}. We now formulate the problems we will be considering. \vspace{2mm} 

We define the slab $\Omega = \{x\in \mathbb{R}^3: 0<x_3<L\}$, with $L>0$, and denote the boundary planes $\Gamma_1=\{ x_3=L\}$ and $\Gamma_2=\{x_3=0\}$. We fix  a frequency $\omega>0$, and assume that the magnetic permeability $\mu$, electric permittivity $\varepsilon$, and conductivity $\sigma$ are {Lipschitz functions} in $\Omega$ such that outside of a compact set, $\mu\equiv \mu_o >0,\varepsilon \equiv \varepsilon_o>0, \sigma \equiv 0$. We set $\gamma=\varepsilon+i\sigma/\omega$ and consider the time-harmonic Maxwell equations on $\Omega$ with boundary data as follows:
\begin{subequations}
 \label{eq:opp-ME}
 \begin{align}
\nabla \wedge E(x) - i \omega\mu(x) H(x) = 0~~\mathrm{in}~\Omega, ~~~& \nu \wedge E\big|_{\Gamma_1} =f, \label{eq:opp-ME1} \\
\nabla \wedge H(x) + i\omega\gamma(x) E(x) = 0~~\mathrm{in}~\Omega,~~~ & \nu \wedge E\big|_{\Gamma_2} =0,  \label{eq:opp-ME2}\\
E\andm H~\mathrm{are~\emph{admissible}~solutions},
\label{eq:opp-ME3}
\end{align}
\end{subequations} 
where $f$ is a compactly supported function in the space $TH^{1/2}_{\mathrm{Div}}(\Gamma_1)$ of tangential vector fields on $\Gamma_1$ 
\begin{eqnarray}
TH^{1/2}_{\mathrm{Div}}(\Gamma_1) = \big\{ F \in H^{1/2}(\Gamma_1)^3 \,|\, \nu\cdot F =0, \nabla_{\partial \Omega}\cdot F \in H^{1/2}(\Gamma_1)\big\}. \label{eq:opp-BCspace}
\end{eqnarray}
The admissibility pertains to guaranteeing uniqueness of the solution on the unbounded domain by prescribing a suitable radiation condition for the fields as $|(x_1,x_2)|\rightarrow \infty$, and this is made precise in Definition \ref{def-admissible} in the appendix, where the well-posedness of Maxwell's equations in this setting is discussed. We show there that the problem \eqref{eq:opp-ME1}-\eqref{eq:opp-ME3} has a unique solution under the conditions stated above, as well as the two following assumptions. \vspace{1mm}

{\it Assumption 1:} We assume that the system \eqref{eq:opp-ME1}-\eqref{eq:opp-ME2} has only the trivial solution if $f=0$.\vspace{1mm}

{\it Assumption 2:} We assume that $k=\omega\sqrt{\mu_o\varepsilon_o}$ is such that $k \neq m\pi/L$ for all $m\in \mathbb{N}$.
\vspace{2mm}

We consider the two cases described above: on the one hand, knowledge of the tangential boundary components $\nu\wedge E$ and $\nu\wedge H$ on (parts of) opposite boundary hyperplanes, and on the other hand, knowledge of these boundary values on the same boundary hyperplane. Consequently, the partial Cauchy data set will be of the form
\begin{equation}
C^D_{\Gamma_2'}(\mu,\varepsilon,\sigma; \omega) = \big\{ \big((\nu\wedge E)|_{\Gamma_1},(\nu\wedge H)|_{\Gamma_2'}\big): (E,H)~\mathrm{solves}~\eqref{eq:opp-ME1}-\eqref{eq:opp-ME3}\big\}, \label{eq:cspd1}
\end{equation}
where $\Gamma_2'\subset \Gamma_2$, in the first case, and in the latter case, with $\Gamma_1'\subset \Gamma_1$,
\begin{equation}
C^S_{\Gamma_1'}(\mu,\varepsilon,\sigma; \omega) = \big\{ \big((\nu\wedge E)|_{\Gamma_1},(\nu\wedge H)|_{\Gamma_1'}\big): (E,H)~\mathrm{solves}~\eqref{eq:opp-ME1}-\eqref{eq:opp-ME3}\big\}.\label{eq:cspd2}
\end{equation}
Our goal is to prove the following theorems. 
\vspace{2mm}

\newdefi{Theorem} \label{thm-same} \emph{Let $\Omega \subset \mathbb{R}^3$ be the slab defined above, and let $\mu_j,\,\varepsilon_j,\,\sigma_j \in C^4(\overline{\Omega}),\,j=1,2$, such that outside of a compact set $B$, $\mu_j\equiv \mu_o >0,\varepsilon_j \equiv \varepsilon_o>0, \sigma_j \equiv 0$. Assume that $\Omega_b:=\Omega\cap B$ is a $C^{1,1}$ domain. Assume furthermore that
\[ \mu_1=\mu_2,~~\gamma_1=\gamma_2~~up~to~order~one~on~\Gamma_1, \]
and that there exist extensions of the parameters to $\mathbb{R}^3$ (of the same name) that belong to $C^4_0(\mathbb{R}^3)$ and are invariant under reflection across the plane $\Gamma_2$.
 Let $\Gamma_1'\subset \Gamma_1$ such that $B\cap\Gamma_1 \subset \Gamma_1'$. Then if for a fixed frequency $\omega>0$ we have $C^S_{\Gamma_1'}(\mu_1,\varepsilon_1,\sigma_1; \omega)=C^S_{\Gamma_1'}(\mu_2,\varepsilon_2,\sigma_2; \omega)$, then $\mu_1=\mu_2$, $\varepsilon_1=\varepsilon_2$, and $\sigma_1=\sigma_2$ in $\Omega$.} \vspace{2mm}
 
\newdefi{Theorem} \label{thm-opp} \emph{Let $\Omega \subset \mathbb{R}^3$ be the slab defined above, and let 
$\mu_j,\,\varepsilon_j,\,\sigma_j \in C^4(\overline{\Omega}),\,j=1,2$, such that outside of a compact set $B$, $\mu_j\equiv \mu_o >0,\varepsilon_j \equiv \varepsilon_o>0, \sigma_j \equiv 0$. Assume that $\Omega_b:=\Omega\cap B$ is a $C^{1,1}$ domain. 
Assume furthermore that
\[ \mu_1=\mu_2,~~\gamma_1=\gamma_2~~up~to~order~one~on~\partial\Omega, \] 
and that there exist extensions of the parameters to $\mathbb{R}^3$ (of the same name) that belong to $C^4_0(\mathbb{R}^3)$ such that $\gamma_1,\mu_1$ are invariant under reflection across the plane $\Gamma_2$, and $\gamma_2,\mu_2$ are invariant under reflection across the plane $\Gamma_1$.
 Let $\Gamma_2'\subset \Gamma_2$ such that $B\cap\Gamma_2 \subset \Gamma_2'$. Then if for a fixed frequency $\omega>0$ we have $C^D_{\Gamma_2'}(\mu_1,\varepsilon_1,\sigma_1; \omega)=C^D_{\Gamma_2'}(\mu_2,\varepsilon_2,\sigma_2; \omega)$, then $\mu_1=\mu_2$, $\varepsilon_1=\varepsilon_2$, and $\sigma_1=\sigma_2$ in $\Omega$.} \vspace{2mm}

This paper is organized as follows. In Section \ref{sec:aug}, we will relate the Maxwell system to a vector Schr\"{o}dinger equation, following an approach originally presented in \cite{OS1996}, with some modifications introduced in \cite{COS2009}. At the end of the section we review the construction of CGO solutions, based on the original method by \cite{SU1987}. 
In Section \ref{sec:same}, we prove Theorem \ref{thm-same}. We first derive an integral formula for pairs of solutions vanishing on the inaccessible hyperplane $\Gamma_2$. We then employ a Runge-type approximation argument to show that it suffices to construct solutions on a bounded domain.  We reflect CGO solutions across $\Gamma_2$ to satisfy the vanishing boundary conditions, and then use these solutions in the integral formula. In the process, products of reflected and non-reflected solutions appear, whose asymptotics need to be studied carefully. The resulting asymptotic expressions will be the same as were obtained in \cite{COS2009}, so that we can refer to that work to finish the uniqueness proof.
Finally, in Section \ref{sec:opp}, we prove Theorem \ref{thm-opp}. We derive a suitable integral formula involving solutions that vanish on different hyperplanes, then construct these solutions by reflecting across $\Gamma_1$ respectively $\Gamma_2$, and use them in the integral formula. In addition to products of reflected and non-reflected solutions, we now also obtain products of solutions reflected across different planes. In order to handle these terms, we adapt an argument used in \cite{KLU2012}, choosing the phase vectors for the CGO solutions suitably so that the products involving reflected solutions exhibit exponential decay. The non-vanishing terms in the limit will again be of the same form as in \cite{COS2009}.

In the appendix, we discuss the necessary conditions to guarantee well-posedness of the direct problem.

\section{Transformation of Maxwell's equations to an elliptic system}\label{sec:aug}
We first modify the Maxwell system, following the approach introduced in \cite{OS1996} and adapted with a slightly different scaling in \cite{COS2009}. Note that \eqref{eq:opp-ME1} implies
\[ \nabla \cdot (\mu H) = \nabla \mu \cdot H + \mu \nabla \cdot H =0 \Rightarrow \nabla \cdot H + \beta \cdot H =0, \]
where $\beta = \nabla \log \mu$, and similarly, with $\alpha=\nabla\log \gamma$,
\[ \nabla \cdot E + \alpha \cdot E =0. \]
We add these two equations to the Maxwell system, and add some suitable terms introducing two scalar potentials $\Phi$ and $\Psi$ to rewrite the system as
\begin{equation}
\big(P(i\nabla)-V\big) u = 0, \label{eq:opp-MEaug}
\end{equation}
with
\[ P
(i\nabla) = i \left(\begin{array}{cccc}
0 &0&0&\nabla\cdot \\
0 & 0 &\nabla & -\nabla \wedge \\
0 & \nabla \cdot& 0 & 0 \\
\nabla & \nabla \wedge & 0 & 0 \end{array} \right), ~~
V = \left(\begin{array}{cccc}
\omega \mu &0 &0&-i \alpha \cdot \\
0& \omega\mu I_3 & -i\alpha  & 0 \\
0 & -i\beta\cdot  & \omega \gamma & 0 \\
-i\beta & 0 &0 & \omega \gamma I_3 \end{array} \right), ~~u = \left(\begin{array}{c} \Phi\\ H \\ \Psi \\E  \end{array} \right).
 \] 
Note that any solution to \eqref{eq:opp-MEaug} that has vanishing first and last components is a solution to Maxwell's equations. In accordance with this block notation, we will from now on write 8-vectors $X=(X_1,X_2,X_3,X_4)$ with $X_1,X_3$ being scalar functions, and $X_2,X_4$ being 3-vectors (corresponding to the magnetic and electric fields, respectively).

For scalars $x$ and $y$ we introduce the notation 
\[ \diag{x,y}= \left(\begin{array}{cccc}
x&0 &0&0\\
0& x I_3 & 0 &0 \\
0 & 0 & y& 0 \\
0 & 0 &0 & y I_3 \end{array} \right). \]
We then define, with $\kappa=\omega(\mu\gamma)^{1/2}$,
\[ W 
= -\frac{i}{2} \left(\begin{array}{cccc}
2i \kappa &0 &0&\alpha \cdot \\
0& 2i \kappa I_3 & \alpha  & \alpha \wedge \\
0 & \beta\cdot  & 2i\kappa & 0 \\
\beta & -\beta\wedge &0 & 2i\kappa I_3 \end{array} \right), \]
and note that if $Y=\diag{\mu^{1/2},\gamma^{1/2}} X$, then
\[ \big(P(i\nabla) -V\big) X=0 \Leftrightarrow \big( P(i\nabla) -W \big) Y=0. \]
The motivation for this rescaling is the following crucial result, which provides a relationship between the augmented Maxwell system and a matrix Schr\"{o}dinger equation.\vspace{2mm}

\newdefi{Lemma} \label{me-schr} \cite[Lemma 1.1]{COS2009} \emph{The operators defined above satisfy
\begin{align*}
\big( \pna - W\big)\big(\pna+W^T\big) &= -\Delta + \tilde{Q}\\
\big( \pna + W^T\big)\big(\pna - W\big) &= -\Delta + \tilde{Q}',
\end{align*}
where $\tilde{Q}$ and $\tilde{Q}'$ are zeroth order matrix potentials, given by
\begin{align*}
\tilde{Q}=&\frac{1}{2}\left(\begin{array}{cccc}
\nabla \cdot \alpha &0&0&0 \\
0 & \nabla \alpha^T+(\nabla \alpha^T)^T -\nabla\cdot \alpha I_3&0 & \\
0 &0 & \nabla\cdot \beta & 0 \\
0& 0 & 0 &\nabla\beta^T+(\nabla\beta^T)^T  -\nabla\cdot\beta I_3 \end{array} \right) \\
&~~~~~ - \left(\begin{array}{cccc}
\kappa -\frac{1}{4}\alpha\cdot\alpha &0&0&-i\nabla \kappa \cdot \\
0 & (\kappa -\frac{1}{4}\alpha\cdot\alpha)I_3 &-i\nabla \kappa & 0 \\
0 & -i\nabla \kappa\cdot  & \kappa -\frac{1}{4}\beta\cdot\beta & 0 \\
-i\nabla \kappa & 0 & 0 & ( \kappa -\frac{1}{4}\beta\cdot\beta)I_3 \end{array} \right), 
\end{align*}
and $\tilde{Q}'$ has the same shape as $\tilde{Q}$ with $\alpha$ and $\beta$ interchanged.\vspace{2mm}
}

For further details concerning the rescaling and properties of the operators we refer to \cite{COS2009}.

\subsection{Review of construction of CGO solutions} \label{sec:cgo}
We briefly summarize the construction of CGO solutions to Maxwell's equations by using the factoring of the Schr\"{o}dinger equation shown in Lemma \ref{me-schr}. The construction follows the classical method developed in \cite{SU1987} and adapted for Maxwell's equations in \cite{OS1996}; see also \cite[Section 2]{COS2009}  for more details.

We define weighted $L^2$ and Sobolev spaces with $-1<\delta<0$ by the following norms:
\begin{align*}
\|f\|_{L^2_\delta}^2 &=\int_{\mathbb{R}^3} \big(1+|x|^2\big)^\delta |f(x)|^2 dx,\\
\|u\|_{H^s_\delta} &=\big\|\big(1+|x|^2\big)^{\delta/2} u\big\|_{H^s(\mathbb{R}^3)}
\end{align*}
In the scalar case, it has been shown that if $q\in L^\infty(\mathbb{R}^3)$ is compactly supported, and $\zeta \in \mathbb{C}^3$ is such that $\zeta\cdot\zeta=0$ and $|\zeta|$ is sufficiently large, then for any $f\in L^2_{\delta+1}(\mathbb{R}^3)$ there is a unique solution $u\in H^2_\delta(\mathbb{R}^3)$ to
\[ (-\Delta -2i\zeta\cdot \nabla +q)u=f, \]
and $u$ satisfies the estimate
\[ \|u\|_{H^s_\delta} \leq C|\zeta|^{s-1}\|f\|_{L^2_{\delta+1}}. \]
We use the analog of this result in the vector case. Recall our extension of the parameters such that $\mu-\mu_o$ and $\gamma-\varepsilon_o$ belong to $ C^4_0(\mathbb{R}^3)$. We set $k=\omega(\mu_o\varepsilon_o)^{1/2}$ and define $Q=k^2+\tilde{Q}$, and note that $Q$ is thus compactly supported. The first operator from Lemma \ref{me-schr} can then be written as 
\begin{equation}
\big( \pna - W\big)\big(\pna+W^T\big) = -(\Delta+k^2)+Q. \label{eq:opp-schr-decomp}
\end{equation}
We first obtain solutions to this second order equation. 
Let $\xi\in \mathbb{R}^3$ be any fixed nonzero vector and let $\zeta\in\CC^3$ be a vector such that $\zeta\cdot \zeta =k^2$ and 
\[ \zeta = \frac{1}{2}\xi + O(\tau), \]
where $\tau\geq 1$ is a parameter controlling the size of $\zeta$.
Let $Z_0=Z_0(\zeta)$ be a vector that is independent of $x$ and bounded with respect to $\tau$. Then we have the following existence result for CGO solutions to the Schr\"{o}dinger equation in $\mathbb{R}^3$ with potential $Q$.\vspace{2mm}

\newdefi{Proposition} \label{cgo} \cite[Proposition 2.1]{COS2009} \emph{Let $-1<\delta<0$, and let $\delta'>0$ such that $-1<\delta+\delta'<0$. There exists a CGO solution to 
\[ \big(-(\Delta+k^2)+Q\big) Z=0~\inm\mathbb{R}^3 \]
that is of the form
\[ Z(x)=e^{i\zeta\cdot x}(Z_0+Z_{-1}(x)+Z_r(x)), \]
such that for $0\leq|\alpha|\leq 2$,
\[ \|\nabla^\alpha Z_{-1}\|_{L^2_\delta}=O(\tau^{-1}),~~\|\nabla^\alpha Z_r\|_{L^2_\delta} = O(\tau^{|\alpha|-(1+\delta')}). \]
This solution satisfies the following asymptotics as $\tau\rightarrow\infty$: if we denote 
\[ \hat{R} = \lim_{\tau\rightarrow\infty} \tau Z_{-1},~~\hat{M}=\lim_{\tau\rightarrow\infty} Z_0, \]
then
\begin{equation}
2i (\hat{\zeta} \cdot \nabla)I_8 \hat{R} = Q \hat{M}. \label{eq:limR}
\end{equation}}
From the factoring \eqref{eq:opp-schr-decomp} it follows that the function
\[ Y= (\pna+W^T)Z \]
then satisfies the first order equation
\[ (\pna-W)Y=0. \]
With
\begin{align}
Y_1&=-P(\zeta) Z_0,\notag\\
Y_0&=-P(\zeta) Z_{-1} + W^T Z_0,\label{eq:opp-ydecomp}\\
Y_r&=(\pna+W^T)Z_{-1}+(\pna -P(\zeta) +W^T) Z_r,\notag
\end{align}
we can write $Y$ as
\[ Y(x)=e^{i\zeta\cdot x} (Y_1+Y_0+Y_r), \]
with the asymptotics
\[ \|Y_1\|_{L^2(U)}=O(\tau),~~\|\nabla^\alpha Y_0\|_{L^2(U)}=O(1),~~\| \nabla^\alpha Y_r\|_{L^2(U)}=O(\tau^{|\alpha|-\delta'}) \]
for all $0\leq \alpha \leq 1$ and bounded open subsets $U$ of $\mathbb{R}^3$.
In order for this function $Y$ to further yield solutions to Maxwell's equations, 
we need $Y_1=Y_3=0$. The next result gives a condition on $Z$ under which this is the case.

\newdefi{Lemma}\label{cgo-me} \cite[Lemma 2.2]{COS2009} \emph{If
\[ \big((-P(\zeta) +k) Z_0\big)_1=\big((-P(\zeta) +k)Z_0\big)_3=0, \]
then if $\tau$ is sufficiently large,
\[ Y_1=Y_3=0. \]}
We will use the following choice of $Z_0$ for $a,b\in\mathbb{C}^3$, which was introduced in \cite{OS1996} and also used in \cite{COS2009}:
\begin{equation}
Z_0= \frac{1}{\tau} \big( \zeta\cdot a, kb, \zeta \cdot b, ka\big).\label{eq:opp-Z0}
\end{equation}

\section{Proof of Theorem \ref{thm-same}} \label{sec:same}

{ The proof of Theorem \ref{thm-same} is divided into several steps. Using the assumption on the Cauchy data sets, we first derive an integral formula in the slab involving the unknown parameters as well as solutions to Maxwell's equations with vanishing boundary conditions on $\Gamma_2$, which is analogous to the integral formula obtained in \cite{COS2009} on a bounded domain. Since the CGO solutions that we want to use grow at infinity, we then need a Runge-type approximation result that allows us to consider the integral formula over a bounded domain. Next, we will reflect the CGO solutions constructed in Section \ref{sec:cgo} across $\Gamma_2$ to achieve the desired boundary conditions. Finally, we plug these solutions into the integral formula and perform the limit $\tau\rightarrow \infty$. The resulting asymptotic expressions are the same as those obtained in \cite{COS2009}, so that we refer to this work to evaluate the limits and obtain partial differential equations for the unknown parameters. A unique continuation result from \cite{COS2009} then shows that the parameters are in fact equal.
}
\subsection{Integral identity} \label{sec:intform-same}

We start by deriving an integral formula for solutions to the augmented Maxwell system. Recall the assumptions of Theorem \ref{thm-same}: we suppose that we can extend the two sets of parameters $\mu_j,\varepsilon_j,\sigma_j,\,j=1,2,$ to $C^4$ functions in all of $\mathbb{R}^3$ in such a fashion that the parameters are constant outside the compact set $B$, so that we have
\[ \mu_j=\mu_o,~~\varepsilon_j=\varepsilon_o,~~\sigma_j=0~~\inm~B^c. \]
Denote $V_j=V(\mu_j,\gamma_j)$, and let $X^1=(0,X^1_2,0,X^1_4)$ be an { admissible} solution to 
\[ \big(P(i\nabla)-V_1\big) X^1 =0~\inm \Omega \]
with $\nu \wedge X^1_4 =0$ on $\Gamma_2$ and $\nu\wedge X^1_4$ compactly supported on $\Gamma_1$. Let $\tilde{X}=(0,\tilde{X}_2,0,\tilde{X}_4)$ solve
\[ \big(P(i\nabla)-V_2\big) \tilde{X}=0~\mathrm{in}~\Omega,~~ \nu \wedge \tilde{X}_4=\nu\wedge X^1_4~\mathrm{on}~\partial \Omega.\]
Then $w=\tilde{X}-X^1$ satisfies $\nu\wedge w_4= 0$ on $\partial \Omega$, and
\begin{equation}
\big(P(i\nabla)-V_2\big) w = (V_2-V_1) X^1~~\mathrm{in}~\Omega. \label{eq:same-w-aux}
\end{equation}
Since $X^1_4$ and $\tilde{X}_4$ have the same tangential boundary values on $\partial \Omega$, it further follows from $C^S_{\Gamma_1'}(\mu_1,\varepsilon_1,\sigma_1; \omega)=C^S_{\Gamma_1'}(\mu_2,\varepsilon_2,\sigma_2; \omega)$ that $\nu\wedge X^1_2=\nu\wedge\tilde{X}_2$ on $\Gamma_1'$, and hence $\nu \wedge w_2  = 0$ on $\Gamma_1'$.
We proceed to show that $w=0$ in $\Omega\backslash B$. Note that since outside $B$ all material parameters are constant, 
\[ V_1=V_2= \left(\begin{array}{cccc}
\omega \mu_o &0 &0&0\\
0 & \omega\mu_o I_3 & 0 & 0 \\
0 & 0 & \omega \varepsilon_o & 0 \\
0 & 0 & 0 & \omega \varepsilon_o I_3 \end{array} \right)\inm \Omega\backslash B, \] 
and it follows from \eqref{eq:same-w-aux} that the nonzero components of $w$ satisfy Maxwell's equations with constant parameters in $\Omega\backslash B$, with partial boundary condition $\nu\wedge w_2=\nu\wedge w_4=0$ on $\Gamma_1'\backslash \ell_1$,
where $\ell_1=B\cap \Gamma_1$. 
Writing out the equations for the components of $w$, we have
\begin{equation}
 \nabla\cdot w_4=0,~~-i\nabla \wedge w_4-\omega\mu_o w_2=0,~~ \nabla \cdot w_2=0,~~i\nabla \wedge w_2-\omega \varepsilon_o w_4=0, \label{eq:same-MEw}
 \end{equation}
and eliminating either of the functions shows that  
we have for $j=1,2$
\[ (-\Delta -k^2) w_j =0 ~ \inm \Omega\backslash B,~~~ \nu\wedge w_j=0 \on \Gamma_1'\backslash\ell_1. \]
Also, $i \omega\mu_o \nu \cdot w_2= \nu \cdot(\nabla \wedge w_4)=-\nabla_{\partial \Omega} \cdot (\nu \wedge w_4) =0$ on $\Gamma_1'\backslash \ell_1$, and similarly for $\nu \cdot w_4$. This implies that $w_2=w_4=0$ on that part of the boundary. Furthermore, given their relationship through Maxwell's equations, the components of $w_4$ are (up to some constants) the normal derivatives of those of $w_2$ on the boundary. So $w_2$ satisfies a homogeneous Helmholtz equation with zero Dirichlet and Neumann boundary condition on a subset of the boundary, and unique continuation now implies $w_2=0$ in $\Omega\backslash B$. By symmetry, the same follows for $w_4$. In particular, we find that $\nu\wedge w_2=\nu\wedge w_4=0$ on $\ell_3=\Omega\cap\partial B$.

Define $\check{V}_2=V(\mu_2,\overline{\gamma_2})$, where $\bar{.}$ denotes the complex conjugate, and 
let $X^2=(0,X^2_2,0,X^2_4)$ solve
\[ \big(P(i\nabla)-\check{V}_2\big) X^2 =0~~~\mathrm{in}~\Omega_b,~~~\nu \wedge X^2_4=0~\mathrm{on}~\ell_2:=B\cap \Gamma_2. \]
Using the integration by parts formula for a bounded domain $D$,
\begin{equation}
\int_{D} \pna U \cdot \bar{u} dx = \int_{\partial D} \pnu U \cdot \bar{u} dS + \int_{D} U \cdot \overline{\pna u} dx, \label{eq:same-ibp}
\end{equation}
as well as the fact that if both $u$ and $U$ have vanishing first and third components, 
\begin{equation}
\int_{\Omega} {V_2} U \cdot \bar{u} dx = \int_{\Omega} U \cdot V_2^T\overline{ u} dx = \int_{\Omega} U\cdot \overline{\check{V}_2 u} dx, \label{eq:same-Vid}
\end{equation}
we compute
\begin{align*}
\int_\Omega (V_2-V_1) X^1 \cdot \overline{X^2} dx &= \int_{\Omega_b} \big(P(i\nabla)-V_2\big)w \cdot \overline{X^2} dx\\
&= \int_{\Omega_b} w\cdot \overline{\big(P(i\nabla) - \check{V}_2 \big) X^2 } dx + \int_{\partial\Omega_b} P(i\nu) w \cdot \overline{X^2} dS\\
&= i \int_{\partial\Omega_b}  -\nu\wedge w_4 \cdot \overline{X^2_2} + \nu\wedge w_2\cdot \overline{X^2_4} dS = i \int_{\ell_2} \nu\wedge w_2 \cdot \overline{X^2_4} dS=0,
\end{align*}
using the boundary values of $w$ and the boundary condition for $X^2$. 
Summarizing, we have shown the following:\vspace{2mm}

\newdefi{Proposition} \emph{Let $X^1=(0,X^1_2,0,X^1_4)$ be an { admissible} solution to 
\[ \big(P(i\nabla)-V_1\big) X^1 =0 \]
in $\Omega$, with $\nu \wedge X^1_4 =0$ on $\Gamma_2$ and $\nu\wedge X^1_4$ compactly supported on $\Gamma_1$. Furthermore, let $X^2=(0,X^2_2,0,X^2_4)$ be a solution to 
\[ \big(P(i\nabla)-\check{V}_2\big) X^2 =0~~~\mathrm{in}~\Omega_b,~~~\nu\wedge X^2_4=0~\mathrm{on}~B\cap \Gamma_2. \]
Then
\begin{equation}
\int_\Omega (V_2-V_1) X^1 \cdot \overline{X^2} dx =0. \label{eq:same-intform}
\end{equation}} \bewende

\subsection{Restricting to a bounded domain - Runge approximation}\label{sec:runge-same}

We introduce the function spaces 
\[ \vcal_j(\Omega_b) = \big\{ u=(0,u_2,0,u_4) \in H^1(\Omega_b)^8:~ \big(P(i\nabla)+\check{V_j}\big)u=0 ~\mathrm{in}~\Omega_b,~~\nu\wedge u_4=0 \on B \cap \Gamma_2 \big\},\]
\[ \wcal_j(\Omega_b) = \left\{ u=(0,u_2,0,u_4) \in H^1(\Omega_b)^8:  
\big(P(i\nabla)-V_j\big)u=0 ~\mathrm{in}~\Omega_b, 
\nu\wedge u_4=0~\mathrm{on}~B\cap\Gamma_2  
\right\},\]
\[ \wcal_j(\Omega) = \left\{ u=(0,u_2,0,u_4) \in H^1(\Omega)^8:~~~\begin{array}{c} \big(P(i\nabla)-V_j\big)u=0 ~\mathrm{in}~\Omega,~~~ 
\nu\wedge u_4=0~\mathrm{on}~\Gamma_2,~\\\nu\wedge u_4 ~\mathrm{compactly~supported~on~} \Gamma_1, \\ $u$\mathrm{~is~admissible~in~the~sense~of~Definition~\ref{def-admissible}}\end{array} \right\}.\]
With this notation, the integral identity holds for $X^1\in \wcal_1(\Omega)$ and $X^2\in \vcal_2(\Omega_b)$. The functions that we will construct to use in the integral formula have exponential growth at infinity, so we want to restrict ourselves to considering solutions on the bounded domain $\Omega_b$ only. To facilitate this, we show the following density result.\vspace{2mm}

\newdefi{Lemma} \label{runge} \emph{ $\wcal_j(\Omega)$ is dense in $\wcal_j(\Omega_b)$ with respect to the $L^2(\Omega_b)$ norm.}\vspace{2mm}

{\it Proof.} Suppose that this density does not hold for $j=1$ (for $j=2$, the argument is analogous). Then the Hahn-Banach Theorem gives the existence of a function $g\in L^2(\Omega)^8$ with $g=0$ in $\Omega\backslash B$ such that for all $u\in \wcal_1(\Omega)$,
\[ \int_\Omega g \cdot \bar{u} \,dx =0, \]
but for some $u_o\in \wcal_1(\Omega_b)$,
\[ \int_\Omega g \cdot \overline{u_o} dx \neq 0. \]
{ We want to replace the function $g$ by $(P(i\nabla)-\check{V}_1)U$, with suitable $U=(0,U_2,0,U_4)$, so that we can integrate by parts. To this end, let $U_2,\,U_4$ be the admissible (in the sense of Definition \ref{def-admissible}) solutions to the following nonhomogeneous Maxwell equations with parameters $\mu_1$ and $\overline{\gamma_1}$ and zero tangential boundary condition for $U_4$,
\[ \nabla\wedge U_4-i \omega \mu_1 U_2 = i g_2,~~~~ \nabla\wedge U_2+i\omega\overline{\gamma_1} U_4=-i g_4~~\mathrm{in}~\Omega,~~~\nu\wedge U_4=0~~\mathrm{on}~\partial \Omega.  \]
Then the second and forth components of $(P(i\nabla)-\check{V}_1)U$ are equal to those of $g$. Since $u_1=u_3=0$, the first and third components of $g$ are not relevant and we can replace $g$ by $(P(i\nabla)-\check{V}_1)U$ 
in the integral. }

We now integrate by parts, using the identity
\begin{equation}
\int_{\Omega} \pna U \cdot \bar{u} dx = \int_{\partial\Omega} \pnu U \cdot \bar{u} dS + \int_{\Omega} U \cdot \overline{\pna u} dx \label{eq:same-Pibp}
\end{equation}
as well as \eqref{eq:same-Vid}.
We thus obtain for all $u\in\wcal_1(\Omega)$
\begin{align*}
0 & = \int_{\Omega} g \cdot \bar{u} \,dx = \int_\Omega \big(P(i\nabla)-\check{V}_1\big)U \cdot  \bar{u}  \, dx 
 = \int_\Omega U \cdot \overline{\big(P(i\nabla)-{V_1}\big) u}\, dx + \int_{\partial \Omega} P(i\nu)U \cdot  \bar{u}  \,dS \\
& = -i \int_{\partial \Omega} \nu\wedge U_4 \cdot  \bar{u}_2  + U_2 \cdot \nu\wedge  \bar{u}_4  \,dS = -i\int_{\Gamma_1} U_2 \cdot \nu\wedge  \bar{u}_4  \,dS = i \int_{\Gamma_1} \nu\wedge U_2 \cdot  \bar{u}_4 \, dS.
\end{align*}
Since $\nu\wedge u_4$ can be an arbitrary smooth function on $\Gamma_1$, we find that $\nu\wedge U_2$ must vanish on $\Gamma_1$.

We proceed to show that $U=0$ in $\Omega\backslash B$. Since $g=0$ and the parameters are constant outside $B$, it follows that $U_2$ and $U_4$ satisfy the homogeneous Maxwell equations in $\Omega\backslash B$. We further have the boundary conditions $\nu\wedge U_2=0$ and $\nu\wedge U_4=0$ on $\Gamma_1\backslash B$. The same unique continuation argument that was used in the derivation of the integral formula for the auxiliary function $w$ applies and yields $U_2=U_4=0$ in $\Omega\backslash B$. In particular, we may conclude that $U=0$ on $\ell_3$. 
This implies for $u_o\in \wcal_1(\Omega_b)$, 
\begin{align*} 
0 &\neq  \int_\Omega g \cdot \overline{u_o} dx =  \int_{\Omega_b} \big(P(i\nabla)-\check{V}_1\big)U \cdot \overline{u_o} dx = \int_{\Omega_b}  U\cdot \overline{\big(P(i\nabla)-{V_1}\big) u_o} dx + \int_{\partial(\Omega_b)} \pnu U \cdot \overline{u_o} dS \\ 
&= \int_{\ell_1\cup\ell_2\cup\ell_3} \nu\wedge U_4 \cdot \overline{u_{o,2}} + U_2\cdot \nu\wedge \overline{u_{o,4}} dS.
\end{align*}
Now the integral over $\ell_3$ vanishes because $U=0$ on $\ell_3$; the first term vanishes on $\ell_1\cup \ell_2$ by the boundary condition for $U$, and the second term vanishes on $\ell_2$ by the boundary condition for $u_o$, and on $\ell_1$, since we saw above that $\nu\wedge U_2=0$ on $\Gamma_1$. Thus, all the boundary terms vanish and we arrive at a contradiction, proving the density of $\wcal_1(\Omega)$ in $\wcal_1(\Omega_b)$. \bewende

\subsection{Constructing CGO solutions that vanish on $\Gamma_2$}
In Section \ref{sec:cgo} we recalled the construction of CGO solutions to Maxwell's equations. The solutions to be used in the integral formula \eqref{eq:same-intform} need to vanish on $\Gamma_2$. This is now achieved by suitably reflecting them across this plane as was also done in \cite{COS2009}. 

We start by picking the complex vectors $\zeta_1$ and $\zeta_2$; our choice is the same as in \cite{COS2009}.
For a fixed vector $\xi=(\xi_1,\xi_2,\xi_3)=(\xi',\xi_3)\in\mathbb{R}^3$ with $|\xi'|>0$, we define the unit vectors $\eta_1$ and $\eta_2$ as 
\[ \eta_1=\frac{1}{|\xi'|} \big(\xi_2,-\xi_1,0\big),~~\eta_2 = \eta_1\wedge \frac{1}{|\xi|}\xi = \frac{1}{|\xi'||\xi|} \big(-\xi_1\xi_3, -\xi_2\xi_3, |\xi'|^2\big). \]
These vectors satisfy $\eta_1\cdot \eta_2=0$ and $\eta_j\cdot \xi=0$ for $j=1,2$. Now we set
\begin{align}
\zeta_1&=\frac{1}{2}\xi + i\Big(\tau^2 + \frac{|\xi|^2}{4}\Big)^{1/2} \eta_1 + \big(\tau^2+k^2\big)\eta_2, \label{eq:same-zeta1} \\
\zeta_2&=-\frac{1}{2}\xi - i\Big(\tau^2 + \frac{|\xi|^2}{4}\Big)^{1/2} \eta_1 + \big(\tau^2+k^2\big)\eta_2,\label{eq:same-zeta2}
\end{align}
where $\tau\geq 1$ is a parameter controlling the size of $|\zeta_j|$. Note that $i\zeta_1+\overline{i\zeta_2}=i\xi$, and $\zeta_j\cdot \zeta_j=k^2$, and as $\tau$ becomes large, we have
\[
\lim_{\tau\rightarrow \infty} \frac{\zeta_1}{\tau}=\lim_{\tau\rightarrow\infty}\frac{\overline{\zeta_2}}{\tau}=\eta_2+i\eta_1 =: \hat{\zeta}.\\
\]
We further set $\check{\zeta}=\overline{\hat{\zeta}}=\eta_2-i\eta_1$. 
With these choices of vectors, let $Z^1,Y^1$ be the CGO solutions for $(\mu_1,\gamma_1)$ with complex phase vector $\zeta_1$ as constructed in Proposition \ref{cgo} and Lemma \ref{cgo-me}, and let $Z^2,Y^2$ be the CGO solutions for $(\mu_2,\overline{\gamma_2})$ with phase $\zeta_2$. { Recall that these are global solutions to the respective equations.}
As in Section \ref{sec:aug} we now denote $X^1=\mathrm{diag}(\mu_1^{-1/2},\gamma_1^{-1/2})Y^1$ and $X_2=\mathrm{diag}(\mu_2^{-1/2},\overline{\gamma}_2^{-1/2})Y^2$, and perform a reflection of these solutions in such a way that the resulting functions also solve Maxwell's equations. 
To this end, we denote the reflection across $\Gamma_2$ in Cartesian coordinates by
\[ x=(x_1,x_2,x_3) \mapsto \dot{x}(x):=(x_1,x_2,-x_3), \]
and let $\dot{\Omega}_b=\{ \dot{x}(x): \,x\in\Omega_b\}$. We also introduce the larger domain
\[ O = \Omega_b \cup \dot{\Omega}_b \cup \mathrm{int}(\Gamma_2\cap\partial \Omega_b). \]
We set 
\[ \dot{I}_4 =\left(\begin{array}{cccc}
-1 &0&0&0\\
0 & 1 &0&0 \\
0 &0& 1 & 0 \\
0 &0 & 0 & -1 \end{array} \right), \]
and define
\[ \dot{X}^j(x)=\diag{\dot{I}_4,-\dot{I}_4} X^j(\dot{x}(x)). \]
It is straightforward to check that these functions satisfy
\[ (\pna - V_1) \dot{X}^1=0,~~~(\pna - \check{V}_2)\dot{X}^2=0 \]
in $\Omega_b$ as well as in $\dot{\Omega}_b$ (recall that by our assumption on the parameters, they are invariant under this reflection), and they satisfy $\nu\wedge \dot{X}^j_4=-\nu\wedge X^j_4$ on $\Gamma_2$ for $j=1,2$. Therefore, $X^j+\dot{X}^j$ is a CGO solution satisfying the required vanishing tangential boundary condition. Summarizing, we have the following existence result.\vspace{2mm}

\newdefi{Proposition} \label{cgo-x-same} \emph{ Given a vector $\xi\in\mathbb{R}^3$ with $|\xi'|>0$, for the sets of parameters $(\mu_1,\gamma_1)$ and $(\mu_2,\overline{\gamma_2})$, there exist CGO solutions $\mathbb{X}^1$ and $\mathbb{X}^2$ satisfying 
\[ (\pna - V_1) \mathbb{X}^1=0,~~~(\pna - \check{V}_2)\mathbb{X}^2=0~~\inm\Omega_b \] 
of the form
\begin{align*}
\mathbb{X}^1&= X^1+\dot{X}^1 = \diag{\mu_1^{-1/2},\gamma_1^{-1/2}}\Big(Y^1(x) + \diag{\ifour,-\ifour} Y^1(\dot{x})\Big), \\
\mathbb{X}^2&= X^2+\dot{X}^2 = \diag{\mu_2^{-1/2},\overline{\gamma}_2^{-1/2}}\Big(Y^2(x) + \diag{\ifour,-\ifour} Y^2(\dot{x})\Big),
\end{align*}
where $Y^j$ are given by \eqref{eq:opp-ydecomp} with the complex vectors $\zeta_j$ defined in \eqref{eq:same-zeta1} and \eqref{eq:same-zeta2} for $j=1,2$, respectively. For $\tau$ large enough, $\mathbb{X}^j$ are solutions to Maxwell's equations in $\Omega_b$, and the tangential components of their electric fields vanish on $\Gamma_2$.}\vspace{2mm}

\subsection{Uniqueness of the parameters}\label{sec:unique-same}

Our next step is to plug the CGO solutions described in Proposition \ref{cgo-x-same} into the integral formula \eqref{eq:same-intform}, and perform the limit $\tau\rightarrow \infty$. 
Recall the shape of the potential $V_j$,
\[V_j = \left(\begin{array}{cccc}
\omega \mu_j &0 &0&-i \alpha_j \cdot \\
& \omega\mu_j I_3 & -i\alpha_j  & 0 \\
0 & -i\beta_j\cdot  & \omega \gamma_j & 0 \\
-i\beta_j & 0 &0 & \omega \gamma_j I_3 \end{array} \right)   
=: \omega \diag{\mu_j,\gamma_j} -i A(\alpha_j,\beta_j). \]
Since for sufficiently large $\tau$, $\XX^j_1=\XX^j_3=0$, we have $A(\alpha_j,\beta_j)\XX^1\cdot\overline{\XX^2}=0$, and therefore
\begin{align*}
(V_2-V_1) \XX^1 \cdot \overline{\XX^2} &= \omega \big( \diag{\mu_2,\gamma_2} - \diag{\mu_1,\gamma_1} \big) \XX^1\cdot \overline{\XX^2} \\
&=\omega\, \diag{\mu_2-\mu_1,\gamma_2-\gamma_1}\diag{\mu_1^{-1/2},\gamma^{-1/2}_1}\big(Y^1+\dot{Y}^1\big) \cdot \diag{\mu_2^{-1/2},\gamma_2^{-1/2}}\overline{\big(Y^2+\dot{Y}^2\big)}\\
&= \omega\,\mathrm{diag}\bigg(\frac{\mu_2-\mu_1}{(\mu_1\mu_2)^{1/2}},\frac{\gamma_2-\gamma_1}{(\gamma_1\gamma_2)^{1/2}}\bigg)\big(Y^1+\dot{Y}^1\big) \cdot \overline{\big(Y^2+\dot{Y}^2\big)},
\end{align*}
where we have again used the notation 
\[ \dot{Y}^j(x)=\diag{\dot{I}_4,-\dot{I}_4} Y^j(\dot{x}(x)). \]
Letting 
\[ \tilde{\mu}=\omega\frac{\mu_2-\mu_1}{(\mu_1\mu_2)^{1/2}},~~~\tilde{\gamma}=\omega\frac{\gamma_2-\gamma_1}{(\gamma_1\gamma_2)^{1/2}}, \]
we use the above identity in the integral formula and obtain the four separate terms 
\begin{equation}
0= \int_{\Omega_b} \td Y^1\cdot \overline{Y^2} + \td Y^1\cdot \overline{\dot{Y}^2} + \td \dot{Y}^1\cdot \overline{Y^2} +\td \dot{Y}^1\cdot \overline{\dot{Y}^2}dx.  \label{eq:same-intformY}
\end{equation}
Recalling that the parameters are invariant under reflection across $\Gamma_2$, the last term can be rewritten as 
\begin{align*}
\int_{\Omega_b}\td \dot{Y}^1\cdot \overline{\dot{Y}^2} dx &= \int_{\Omega_b} \td \diag{\ifour,-\ifour} Y^1(\dot{x})\cdot \diag{\ifour,-\ifour} \overline{Y^2(\dot{x})} dx \\
&=\int_{\Omega_b} \diag{\tilde{\mu}(\dot{x}),\tilde{\gamma}(\dot{x})} Y^1(\dot{x})\cdot \overline{Y^2(\dot{x})}dx = \int_{\dot{\Omega}_b} \td Y^1\cdot \overline{Y^2} dx,
\end{align*}
and similarly we get for the third term in \eqref{eq:same-intformY},
\begin{align*}
\int_{\Omega_b}\td \dot{Y}^1\cdot \overline{Y^2} dx &= \int_{\Omega_b} \td \diag{\ifour,-\ifour} Y^1(\dot{x})\cdot \overline{Y^2(x)} dx \\
&=\int_{\Omega_b}\diag{\tilde{\mu}(\dot{x}),\tilde{\gamma}(\dot{x})} Y^1(\dot{x})\cdot \overline{\diag{\ifour,-\ifour} Y^2}dx = \int_{\dot{\Omega}_b} \td Y^1\cdot \overline{\dot{Y}^2} dx,
\end{align*}
so that \eqref{eq:same-intformY} becomes (recall that $O=\Omega_b\cup\dot{\Omega}_b\cup(\mathrm{int}\Gamma_2\cap \partial\Omega_b)$)
\begin{equation}
0= \int_O \td Y^1\cdot \overline{(Y^2+\dot{Y}^2)}  dx  . \label{eq:same-intformO}
\end{equation}
Now we substitute $Y^1=\big(\pna + W^T_1\big)Z^1$. Using the identity
\[ \pna \tdg Z^1 = \td \pna Z^1 + \pna(\tdg) Z^1, \]
we can write 
\begin{equation}
\td Y^1 = \td \big(\pna+W_1^T\big) Z^1 = \big( \pna \tdg- \pna[\tdg] +\td W_1^T\big) Z^1. \label{eq:same-tdy1id}
\end{equation}
An integration by parts yields 
\begin{align*}
 \int_O \pna \tdg Z^1 \cdot \overline{(Y^2+\dot{Y}^2)} dx =& \int_O \tdg Z^1 \cdot \overline{\pna(Y^2+\dot{Y}^2)} dx \\
&~~~~~~~~~~+\int_{\partial O} \pnu \tdg Z^1 \cdot \overline{(Y^2+\dot{Y}^2)} dS.  
\end{align*}
By assumption, the material parameters are equal up to order one on $\Gamma_1$, and equal to $\mu_o,\varepsilon_o$ on $\partial B$, so that $\tilde{\mu}=\tilde{\gamma}=0$ on $\partial O$, hence the boundary integral vanishes. 
In the volume integral on the right-hand side, we use the fact that $(\pna - \check{W}_2)(Y^2+\dot{Y}^2)=0$ to obtain
\[  \int_O \pna \tdg Z^1 \cdot \overline{(Y^2+\dot{Y}^2)} dx = \int_O \tdg Z^1 \cdot \overline{\check{W}_2(Y^2+\dot{Y}^2)} dx. \]
Thus, substituting \eqref{eq:same-tdy1id} in \eqref{eq:same-intformO} and using this identity yields
\begin{align*}
0&= \int_O \tdg Z^1 \cdot \overline{\check{W}_2(Y^2+\dot{Y}^2)} - \pna(\tdg) Z^1 \cdot \overline{(Y^2+\dot{Y}^2)} + \td W_1^TZ^1   \cdot \overline{(Y^2+\dot{Y}^2)}dx\\
&= \int_O \Big( \overline{\check{W}_2^T}\tdg + \td W_1^T- \pna(\tdg) \Big)Z^1 \cdot \overline{(Y^2+\dot{Y}^2)} dx.
\end{align*}
Setting $\hat{\mu}=\omega\frac{\mu_1+\mu_2}{(\mu_1\mu_2)^{1/2}}$ and $\hat{\gamma}=\omega\frac{\gamma_1+\gamma_2}{(\gamma_1\gamma_2)^{1/2}}$, we write the operator multiplying $Z^1$ as 
\begin{equation}
\kappa_2 \tdg + \kappa_1 \td + i  \left(\begin{array}{cccc}
0 &0&0&\nabla \hat{\gamma}\cdot\\
0 & 0 &\nabla \hat{\gamma} &\nabla \hat{\gamma}\wedge \\
0 &\nabla \hat{\mu}\cdot & 0 & 0 \\
\nabla \hat{\mu}&-\nabla \hat{\mu}\wedge & 0 & 0 \end{array} \right) - i  \left(\begin{array}{cccc}
0 &0&0&\nabla \tilde{\gamma}\cdot\\
0 & 0 &\nabla \tilde{\gamma} &-\nabla \tilde{\gamma}\wedge \\
0 &\nabla \tilde{\mu}\cdot & 0 & 0 \\
\nabla \tilde{\mu}&\nabla \tilde{\mu}\wedge & 0 & 0 \end{array} \right) =:U, \label{eq:U}
\end{equation}
and can thus write the resulting integral identity concisely as
\begin{equation}
\int_O U Z^1 \cdot \overline{(Y^2+\dot{Y}^2)} dx = 0. \label{eq:same-intformU}
\end{equation}
Written separately for the terms involving $Y^2$ and those involving $\dot{Y}^2$, the above computations give the identities
\[ \int_O \td Y^1\cdot \overline{Y^2}dx = \int_O U Z^1\cdot \overline{Y^2} dx,~~~\int_O \td Y^1\cdot \overline{\dot{Y}^2}dx = \int_O U Z^1\cdot \overline{\dot{Y}^2} dx. \]
Now we recall the asymptotics for $Z^1$ and $Y^2$ in Proposition \ref{cgo} and the discussion following it. Writing
\[ Z^1=e^{i\zeta_1\cdot x}(Z^1_0+Z^1_{-1}+Z^1_r),~~Y^2=e^{i\zeta_2\cdot x}(Y^2_1+Y^2_0+Y^2_r), \]
we have 
\begin{gather*}
\| Z^1_{-1}\|_{L^2_\delta}=O(\tau^{-1}),~~\|Z^1_r\|_{L^2_\delta} = O(\tau^{-(1+\delta')}),\\
 \|Y^2_1\|_{L^2(O)}=O(\tau),~~\| Y^2_0\|_{L^2(O)}=O(1),~~\| Y^2_r\|_{L^2(O)}=O(\tau^{-\delta'}), 
 \end{gather*}
for sufficiently small $\delta'>0$, so that taking the limit $\tau\rightarrow\infty$ in \eqref{eq:same-intformU} and recalling $i\zeta_1+\overline{i\zeta_2}=i\xi$, we obtain
\begin{align}
0&= \lim_{\tau\rightarrow\infty}\int_O Ue^{i\zeta_1\cdot x}(Z^1_0+Z^1_{-1}+Z^1_r) \cdot \overline{e^{i\zeta_2\cdot x}(Y^2_1+Y^2_0+Y^2_r)} + U Z^1\cdot \overline{\dot{Y}^2}dx \notag\\
&= \lim_{\tau\rightarrow\infty} \int_O e^{i\xi\cdot x} U Z^1_0\cdot \overline{Y^2_1} dx \label{eq:same-lim1} \\
&~~~~+ \lim_{\tau\rightarrow\infty} \int_O e^{i\xi\cdot x} U Z^1_0\cdot \overline{Y^2_0} dx\label{eq:same-lim2}  \\
&~~~~+ \lim_{\tau\rightarrow\infty} \int_O e^{i\xi\cdot x} U Z^1_{-1}\cdot \overline{Y^2_1} dx\label{eq:same-lim3} \\
&~~~~+ \lim_{\tau\rightarrow\infty} \int_O U Z^1\cdot \overline{\dot{Y}^2} dx.\label{eq:same-lim4} 
\end{align}
The same expression was obtained in \cite{COS2009}, and each of these limits was computed there. We summarize the results below.\vspace{2mm}

\newdefi{Lemma} \label{same-limlem} {\it Recall the choices of $Z^j_0$ and $Y^j$ in \eqref{eq:opp-Z0} and \eqref{eq:opp-ydecomp}. The limit in \eqref{eq:same-lim1} is
\begin{align}
\eqref{eq:same-lim1} =& -k \big( (\hat{\zeta}\cdot \overline{b_2})(\hat{\zeta}\cdot b_1)+(\hat{\zeta}\cdot \overline{a_2})(\hat{\zeta}\cdot a_1)\big) \int_O e^{i\xi\cdot x} (\kappa_1+\kappa_2)(\tilde{\mu}+\tilde{\gamma}) dx \notag \\
 &~~~+ \omega (\hat{\zeta}\cdot \overline{b_2})(\hat{\zeta}\cdot b_1) \int_O e^{i\xi\cdot x}(-\Delta)\Big(\frac{\mu_1}{\mu_2}\Big)^{1/2} dx + \omega (\hat{\zeta}\cdot \overline{a_2})(\hat{\zeta}\cdot a_1) \int_O e^{i\xi\cdot x}(-\Delta)\Big(\frac{\gamma_1}{\gamma_2}\Big)^{1/2} dx\notag \\
 &~~~+ 2 k\omega \int_O e^{i\xi\cdot x}\bigg(i (\hat{\zeta}\cdot \overline{b_2})(\hat{\zeta}\wedge a_1)\cdot \nabla\Big(\frac{\mu_2}{\mu_1}\Big)^{1/2} - i (\hat{\zeta}\cdot \overline{a_2})(\hat{\zeta}\wedge b_1)\cdot \nabla\Big(\frac{\gamma_2}{\gamma_1}\Big)^{1/2} \bigg) dx \notag \\
&~~~+ 2 k\omega \int_O e^{i\xi\cdot x}\bigg(i (\hat{\zeta}\cdot {b_1})(\hat{\zeta}\wedge \overline{a_2})\cdot \nabla\Big(\frac{\mu_1}{\mu_2}\Big)^{1/2} - i (\hat{\zeta}\cdot {a_1})(\hat{\zeta}\wedge \overline{b_2})\cdot \nabla\Big(\frac{\gamma_1}{\gamma_2}\Big)^{1/2} \bigg) dx.\label{eq:same-limm1}
\end{align}
With 
\[ \check{R}_2 = ( \check{R}_{2,1}, \check{R}_{2,2}, \check{R}_{2,3}, \check{R}_{2,4}) = \lim_{\tau\rightarrow\infty} \tau Z^2_{-1},~~~\hat{R}_2=\overline{\check{R}_2}, \]
where the first and third components are scalars and the second and fourth components are 3-vectors, the limit in \eqref{eq:same-lim2} is
\begin{align}
\eqref{eq:same-lim2} =& \int_O\exi \kappa_2\bigg( (\hat{\zeta}\cdot \overline{a_2})(\hat{\zeta}\cdot a_1)(\kappa_2\tilde{\gamma}+\kappa_1\tilde{\mu}) +(\hat{\zeta}\cdot \overline{b_2})(\hat{\zeta}\cdot b_1)(\kappa_2\tilde{\mu}+\kappa_1\tilde{\gamma})\bigg) dx \notag\\
 &~~~- \omega \int_O\exi \bigg( (\hat{\zeta}\cdot \overline{a_2})(\hat{\zeta}\cdot a_1)\frac{\nabla \gamma_2}{\gamma_2} \cdot \nabla  \Big(\frac{\gamma_1}{\gamma_2}\Big)^{1/2}+(\hat{\zeta}\cdot \overline{b_2})(\hat{\zeta}\cdot b_1) \frac{\nabla\mu_2}{\mu_2}\cdot \nabla \Big(\frac{\mu_1}{\mu_2}\Big)^{1/2}\bigg) dx \notag \\
 &~~~ - \int_O\exi \Big( (\hat{\zeta}\cdot \hat{R}_{2,4})(\hat{\zeta}\cdot a_1) (\kappa_2\tilde{\gamma}+\kappa_1\tilde{\mu}) +  (\hat{\zeta}\cdot \hat{R}_{2,2})(\hat{\zeta}\cdot b_1) (\kappa_2\tilde{\mu}+\kappa_1\tilde{\gamma}) \Big) dx \notag\\
&~~~-2i\omega \int_O\exi \bigg( (\hat{\zeta}\cdot {b_1})(\hat{R}_{2,4}\wedge \hat{\zeta}) \cdot \nabla  \Big(\frac{\mu_1}{\mu_2}\Big)^{1/2}-(\hat{\zeta}\cdot {a_1})(\hat{R}_{2,2}\wedge\hat{\zeta}) \cdot \nabla \Big(\frac{\gamma_1}{\gamma_2}\Big)^{1/2}\bigg) dx \notag \\
&~~~ +2i\omega \int_O \exi \bigg( (\hat{\zeta}\cdot b_1)(\hat{\zeta}\cdot \nabla) \hat{R}_{2,3} \Big[ \Big(\frac{\mu_1}{\mu_2}\Big)^{1/2} - 1\Big] + (\hat{\zeta}\cdot a_1)(\hat{\zeta}\cdot \nabla)\hat{R}_{2,1}\Big[\Big(\frac{\gamma_1}{\gamma_2}\Big)^{1/2}-1\Big]\bigg) dx.\label{eq:same-limm2}
\end{align}
With
\[ \hat{R}_1 = ( \hat{R}_{1,1}, \hat{R}_{1,2}, \hat{R}_{1,3}, \hat{R}_{1,4}) = \lim_{\tau\rightarrow\infty} \tau Z^1_{-1},\]
the limit in \eqref{eq:same-lim3} is 
\begin{align}
\eqref{eq:same-lim3} =& - \int_O\exi \bigg( (\hat{\zeta}\cdot \overline{b_2})(\hat{\zeta}\cdot \hat{R}_{1,2})(\kappa_2\tilde{\gamma}+\kappa_1\tilde{\mu}) +(\hat{\zeta}\cdot \overline{a_2})(\hat{\zeta}\cdot \hat{R}_{1,4})(\kappa_2\tilde{\mu}+\kappa_1\tilde{\gamma})\bigg) dx \notag\\
&~~~ -2i\omega \int_O\exi \bigg( (\hat{\zeta}\cdot \overline{a_2})( \hat{\zeta}\wedge\hat{R}_{1,2}) \cdot \nabla  \Big(\frac{\gamma_2}{\gamma_1}\Big)^{1/2}-(\hat{\zeta}\cdot \overline{b_2})(\hat{\zeta}\wedge \hat{R}_{1,4}) \cdot \nabla \Big(\frac{\mu_2}{\mu_1}\Big)^{1/2}\bigg) dx \notag \\
&~~~ +2i\omega \int_O \exi \bigg( (\hat{\zeta}\cdot \overline{b_2})(\hat{\zeta}\cdot \nabla) \hat{R}_{1,3} \Big[ \Big(\frac{\mu_1}{\mu_2}\Big)^{1/2} - 1\Big] + (\hat{\zeta}\cdot \overline{a_2})(\hat{\zeta}\cdot \nabla)\hat{R}_{1,1}\Big[\Big(\frac{\gamma_1}{\gamma_2}\Big)^{1/2}-1\Big]\bigg) dx.\label{eq:same-limm3}
\end{align}
The limit in \eqref{eq:same-lim4} is
\begin{equation}
\lim_{\tau\rightarrow\infty} \int_O U Z^1\cdot \overline{\dot{Y}^2} dx =0. \label{eq:same-limm4}
\end{equation}
\bewende
}
Using proper choices of the vectors $a_1,b_1$ and $a_2,b_2$ determining $Z^1_0$ and $Z^2_0$, and combining \eqref{eq:same-limm1}-\eqref{eq:same-limm4} as well as the identity \eqref{eq:limR} which both $\hat{R}_1$ and $\check{R}_2$ satisfy,  
we arrive at a set of differential equations involving the unknown parameters:\vspace{2mm}

\newdefi{Proposition} \label{pdes} \cite[Proposition 5.1]{COS2009} {\it Let $\xi\in\RR^3$. If we let $b_1=\overline{b_2}=\check{\zeta}$ and $a_1=\overline{a_2}=\hat{\zeta}$, then 
\[
\int_O\exi \bigg( \frac{1}{2} \nabla \cdot (\beta_2-\beta_1) +\frac{1}{4}\big( \beta_2\cdot\beta_2-\beta_1\cdot \beta_1\big)+ \kappa_1^2-\kappa_2^2\bigg) dx=0.
\]
If we let $a_1=\overline{a_2}=\check{\zeta}$ and $b_1=\overline{b_2}=\hat{\zeta}$, then 
\[
\int_O\exi \bigg( \frac{1}{2} \nabla \cdot (\alpha_2-\alpha_1) +\frac{1}{4}\big( \alpha_2\cdot \alpha_2-\alpha_1\cdot \alpha_1\big)+ \kappa_1^2-\kappa_2^2\bigg) dx=0.
\]\bewende}
Thus, we have the following set of equations in terms of the parameters $\mu_1,\gamma_1,\mu_2,\gamma_2$:
\begin{align*}
-\frac{1}{2} \Delta (\log \mu_1-\log\mu_2) -\frac{1}{4}\nabla(\log\mu_1+\log\mu_2 )\cdot\nabla(\log\mu_1-\log\mu_2)+ \omega^2(\mu_1\gamma_1-\mu_2\gamma_2) &=0, \\
-\frac{1}{2} \Delta (\log \gamma_1-\log\gamma_2) -\frac{1}{4}\nabla(\log\gamma_1+\log\gamma_2 )\cdot\nabla(\log\gamma_1-\log\gamma_2)+ \omega^2(\mu_1\gamma_1-\mu_2\gamma_2) &=0.
\end{align*}
Setting $u=(\gamma_1/\gamma_2)^{1/2}$ and $v=(\mu_1/\mu_2)^{1/2}$, we can rewrite these equations as 
\begin{align*}
- \Delta (\log v) -(\mu_1\mu_1)^{-1/2}\nabla(\mu_1\mu_1)^{1/2}\cdot \nabla(\log v)+ \omega^2(\mu_1\gamma_1-\mu_2\gamma_2) &=0, \\
-\Delta (\log u) -(\gamma_1\gamma_2)^{-1/2}\nabla(\gamma_1\gamma_2 )^{1/2}\cdot\nabla(\log u)+ \omega^2(\mu_1\gamma_1-\mu_2\gamma_2) &=0. 
\end{align*}
After multiplying the first equation by $(\mu_1\mu_2)^{1/2}$ and the second one by $(\gamma_1\gamma_2)^{1/2}$ and combining terms using the product rule, we arrive at 
\begin{align*}
- \nabla\cdot (\mu_2\nabla v)+ \omega^2\mu_2^2\gamma_2(u^2v^2-1)v&=0,\\
-\nabla \cdot (\gamma_2\nabla u) + \omega^2\mu_2\gamma_2^2(u^2v^2-1)u&=0.
\end{align*}
By our assumptions on the parameters, we further have the following boundary conditions for $u$ and $v$:
\[ u = v = 1,~~ \frac{\partial u}{\partial \nu} = \frac{\partial v}{\partial \nu} =0~~\on \partial O. \]
The uniqueness proof is now completed by using the following unique continuation property that was proved in \cite[Lemma 5.2]{COS2009}:\vspace{2mm}

\newdefi{Proposition} \label{uniquecont} {\it Let $a,b\in C^2(\overline{O})$ be non-vanishing complex valued functions with positive real parts, and let $p,q\in L^\infty(O)$ be complex valued functions. Suppose that the functions $u,v \in C^2(\overline{O})$ satisfy
\[
\left. \begin{array}{r}-\nabla\cdot (a \nabla u) + p(u^2v^2-1)u=0\\-\nabla\cdot(b\nabla v) +q(u^2v^2-1)v=0 \end{array}\right\} 
~~\inm~O,~~~~~\left.\begin{array}{r} u=v=1 \\  \frac{\partial u}{\partial \nu} = \frac{\partial v}{\partial \nu} =0 \end{array}\right\}~~\on~\partial O.~~~~~~~~~~~~~~~
\]
Then $u\equiv 1\equiv v$ in $O$. \bewende}

\section{Proof of Theorem \ref{thm-opp}}\label{sec:opp}

The proof of Theorem \ref{thm-opp}, in which data are available on different boundary hyperplanes, is similar to that of Theorem \ref{thm-same}, with some notable differences. Most importantly, since we have data on different planes, in the derivation of a suitable integral formula we now need one solution to vanish on $\Gamma_1$ and the other to vanish on $\Gamma_2$, in order to eliminate all boundary terms. To construct these solutions, we need to employ reflections across both planes, leading to more complicated mixed terms when plugging the solutions into the integral formula. In order to handle these terms, we employ a different choice of vectors $\zeta_1$ and $\zeta_2$ than was done above. Namely, the vectors chosen here will result in exponentially decaying terms as $\tau\rightarrow \infty$. This idea was also used in \cite{KLU2012} in studying the inverse problem for a magnetic Schr\"{o}dinger operator in a slab. The nonvanishing terms when taking the limit will be the same as in the proof of Theorem \ref{thm-same}, so that the final argument is identical.

\subsection{Integral identity} \label{sec:intform-opp}
Analogously to the integral formula \eqref{eq:same-intform}, we can prove the following identity. \vspace{2mm}

\noindent\newdefi{Proposition} \emph{Let $X^1=(0,X^1_2,0,X^1_4)$ be an { admissible} solution to 
\[ \big(P(i\nabla)-V_1\big) X^1 =0 \]
in $\Omega$, with $\nu \wedge X^1_4 =0$ on $\Gamma_2$ and $\nu\wedge X^1_4$ compactly supported on $\Gamma_1$. Furthermore, let $X^2=(0,X^2_2,0,X^2_4)$ be a solution to 
\[ \big(P(i\nabla)-\check{V}_2\big) X^2 =0~~~\mathrm{in}~\Omega_b,~~~\nu\wedge X^2_4=0~\mathrm{on}~\ell_1:=B\cap \Gamma_1. \]
Then
\begin{equation}
\int_\Omega (V_2-V_1) X^1 \cdot \overline{X^2} dx =0. \label{eq:opp-intform}
\end{equation}} 

\noindent{\it Proof.} Let $X^1=(0,X^1_2,0,X^1_4)$ be as in the statement of the proposition, and let $\tilde{X}=(0,\tilde{X}_2,0,\tilde{X}_4)$ { be an admissible solution to}
\[ \big(P(i\nabla)-V_2\big) \tilde{X}=0~\mathrm{in}~\Omega,~~ \nu \wedge \tilde{X}_4=\nu\wedge X^1_4~\mathrm{on}~\partial \Omega.\]
Then $w=\tilde{X}-X^1$ satisfies $\nu\wedge w_4= 0$ on $\partial \Omega$, and
\begin{equation}
\big(P(i\nabla)-V_2\big) w = (V_2-V_1) X^1~~\mathrm{in}~\Omega. \label{eq:opp-w-aux}
\end{equation}
Since $X^1_4$ and $\tilde{X}_4$ have the same tangential boundary values on $\partial \Omega$, it follows from $C^D_{\Gamma_2'}(\mu_1,\varepsilon_1,\sigma_1; \omega)=C^D_{\Gamma_2'}(\mu_2,\varepsilon_2,\sigma_2; \omega)$ that $\nu\wedge X^1_2=\nu\wedge\tilde{X}_2$ on $\Gamma_2'$, and hence $\nu \wedge w_2 = \nu \wedge w_4 = 0$ on $\Gamma_2'$.
Analogously to the argument in Section \ref{sec:intform-same}, we find that $w=0$ in $\Omega\backslash B$. 
Taking $X^2$ as in the statement of the proposition and using the identities \eqref{eq:same-ibp} and \eqref{eq:same-Vid}, we compute
\begin{align*}
\int_\Omega (V_2-V_1) X^1 \cdot \overline{X^2} dx &= \int_{\Omega_b} \big(P(i\nabla)-V_2\big)w \cdot \overline{X^2} dx\\
&= \int_{\Omega_b} w\cdot \overline{\big(P(i\nabla) - \check{V}_2 \big) X^2 } dx + \int_{\partial\Omega_b} P(i\nu) w \cdot \overline{X^2} dS\\
&= i \int_{\partial\Omega_b}  -\nu\wedge w_4 \cdot \overline{X^2_2} + \nu\wedge w_2\cdot \overline{X^2_4} dS = i \int_{\ell_1} \nu\wedge w_2 \cdot \overline{X^2_4} dS=0,
\end{align*}
using the boundary values of $w$ and the boundary condition for $X^2$. \bewende

Recalling the definition of the function spaces in Section \ref{sec:runge-same}, this integral formula holds for $X^1\in \wcal_1(\Omega)$ and $X^2\in \vcal_2^{{D}}(\Omega_b)$, where
\[ \vcal_2^{{D}}(\Omega_b) = \big\{ u=(0,u_2,0,u_4) \in H^1(\Omega_b)^8:~ \big(P(i\nabla)+\check{V_j}\big)u=0 ~\mathrm{in}~\Omega_b,~~\nu\wedge u_4=0 \on B \cap \Gamma_1 \big\}.\]
By Lemma \ref{runge}, it also holds for $X^1\in \wcal_1(\Omega_b)$.

\subsection{Construction of CGO solutions with vanishing boundary values}

We now need the CGO solutions to have vanishing tangential boundary values on $\Gamma_2$ respectively $\Gamma_1$, and we achieve this by suitably reflecting each solution across the respective plane. We start by choosing the complex vectors $\zeta_1$ and $\zeta_2$. The choice now is different from that in the proof of Theorem \ref{thm-same}, so as to be able to control the products of solutions that were reflected across different planes.
For a fixed vector $\xi\in\mathbb{R}^3$ with $|\xi'|>0$, we define the unit vectors $\eta_1$ and $\eta_2$ as before as 
\[ \eta_1=\frac{1}{|\xi'|} \big(\xi_2,-\xi_1,0\big),~~\eta_2 = \eta_1\wedge \frac{1}{|\xi|}\xi = \frac{1}{|\xi'||\xi|} \big(-\xi_1\xi_3, -\xi_2\xi_3, |\xi'|^2\big). \]
Now we set 
\begin{align}
\zeta_1&=~\frac{1}{2}\xi + \Big(\tau^2 - \frac{|\xi|^2}{4}\Big)^{1/2} \eta_1 -i \big(\tau^2-k^2\big)\eta_2, \label{eq:opp-zeta1} \\
\zeta_2&=-\frac{1}{2}\xi + \Big(\tau^2 - \frac{|\xi|^2}{4}\Big)^{1/2} \eta_1 + i\big(\tau^2-k^2\big)\eta_2,\label{eq:opp-zeta2}
\end{align}
with $\tau\geq 1$. The same choice of vectors was used in \cite{KLU2012}, where the inverse problem for a magnetic Schr\"{o}dinger operator in the slab was studied. Note that as before we have $i\zeta_1+\bar{i\zeta_2}=i\xi$, and $\zeta_j\cdot \zeta_j=k^2$, and as $\tau$ becomes large,
\[
\lim_{\tau\rightarrow \infty} \frac{\zeta_1}{\tau}=\lim_{\tau\rightarrow\infty}\frac{\overline{\zeta_2}}{\tau}=\eta_1-i\eta_2 =: \tilde{\zeta}.\\
\]
With these choices of vectors, let $Z^1,Y^1$ be the CGO solutions for parameters $(\mu_1,\gamma_1)$ with complex phase vector $\zeta_1$ as constructed in Proposition \ref{cgo} and Lemma \ref{cgo-me}, and let $Z^2,Y^2$ be the CGO solutions for $(\mu_2,
\overline{\gamma_2})$ with phase $\zeta_2$. 
We denote $X^1=\diag{\mu_1^{-1/2},\gamma_1^{-1/2}}Y^1$ and $X^2=\diag{\mu_2^{-1/2},\overline{\gamma}_2^{-1/2}}Y^2$, and reflect these functions across $\Gamma_2$ respectively $\Gamma_1$. 
Recall that we denote the reflection across $\Gamma_2$ in Cartesian coordinates by
\[ x=(x_1,x_2,x_3) \mapsto \dot{x}(x):=(x_1,x_2,-x_3). \]
Similarly, the reflection across $\Gamma_1$ in Cartesian coordinates is denoted by
\[ x=(x_1,x_2,x_3) \mapsto \ddot{x}(x):=(x_1,x_2,2L-x_3), \]
and we set $\ddot{\Omega}_b=\{ \dot{x}(x): \,x\in\Omega_b\}$. 
With
\[ \dot{I}_4 =\left(\begin{array}{cccc}
-1 &0&0&0\\
0 & 1 &0&0 \\
0 &0& 1 & 0 \\
0 &0 & 0 & -1 \end{array} \right) \]
as before, we define
\[ \dot{X}^1(x)=\diag{\dot{I}_4,-\dot{I}_4} X^1(\dot{x}(x)),~~~\ddot{X}^2(x)=\diag{\dot{I}_4,-\dot{I}_4} X^2(\ddot{x}(x)). \]
One can again check that these functions satisfy
\[ (\pna - V_1) \dot{X}^1=0,~~~(\pna - \check{V}_2)\ddot{X}^2=0 \]
in $\Omega_b$, as well as in $\dot{\Omega}_b$ respectively $\ddot{\Omega}_b$ (recall that we extended the two sets of parameters suitably across the planes), and they satisfy 
\[
\nu\wedge \dot{X}^1_4=-\nu\wedge X^1_4\on \Gamma_2,~\andm~ \nu\wedge \ddot{X}^2_4=-\nu\wedge X^2_4\on\Gamma_1.
\]
Therefore, $X^1+\dot{X}^1$ and $X^2+\ddot{X}^2$ are CGO solutions satisfying the required boundary conditions for use in the integral formula \eqref{eq:opp-intform}. We summarize the construction below.\vspace{2mm}

\newdefi{Proposition} \label{cgo-x-opp} \emph{ Given a vector $\xi\in\mathbb{R}^3$ with $|\xi'|>0$, for the sets of parameters $(\mu_1,\gamma_1)$ and $(\mu_2,\overline{\gamma_2})$, there exist CGO solutions $\mathbb{X}^1$ and $\mathbb{X}^2$ satisfying 
\[ (\pna - V_1) \mathbb{X}^1=0,~~~(\pna - \check{V}_2)\mathbb{X}^2=0 \]
in $\Omega_b$ of the form
\begin{align*}
\mathbb{X}^1&= X^1+\dot{X}^1 = \diag{\mu_1^{-1/2},\gamma_1^{-1/2}}\Big(Y^1(x) + \diag{\ifour,-\ifour} Y^1(\dot{x})\Big), \\
\mathbb{X}^2&= X^2+\ddot{X}^2 = \diag{\mu_2^{-1/2},\overline{\gamma}_2^{-1/2}}\Big(Y^2(x) + \diag{\ifour,-\ifour} Y^2(\ddot{x})\Big),
\end{align*}
where $Y^j$ are given by \eqref{eq:opp-ydecomp} with the complex vectors $\zeta_j$ defined in \eqref{eq:opp-zeta1} and \eqref{eq:opp-zeta2} for $j=1,2$, respectively. For $\tau$ large enough, $\mathbb{X}^1$ and $\mathbb{X}^2$ are solutions to Maxwell's equations in $\Omega_b$, and the tangential components of their electric fields vanish on $\Gamma_2$ and $\Gamma_1$, respectively.}\vspace{2mm}

\subsection{Uniqueness of the parameters}

Our next step is to plug the CGO solutions described in Proposition \ref{cgo-x-opp} into the integral formula \eqref{eq:opp-intform}, and perform the limit $\tau\rightarrow \infty$. As in Section \ref{sec:unique-same}, we can write 
\[ (V_2-V_1) \XX^1 \cdot \overline{\XX^2} = \diag{\tilde{\mu},\tilde{\gamma}} \big(Y^1+\dot{Y}^1\big) \cdot \overline{\big(Y^2+\ddot{Y}^2\big)}, \]
and using this in the integral formula \eqref{eq:opp-intform}, we obtain
\begin{equation}
0= \int_{\Omega_b} \td Y^1\cdot \overline{Y^2} + \td Y^1\cdot \overline{\ddot{Y}^2} + \td \dot{Y}^1\cdot \overline{Y^2} +\td \dot{Y}^1\cdot \overline{\ddot{Y}^2}dx.  \label{eq:opp-intformY}
\end{equation}
We now compute the limit of each of these terms as $\tau \rightarrow \infty$. We first rewrite the first two terms, substituting $Y^1=\big(\pna + W^T_1\big)Z^1$, as was also done in Section \ref{sec:unique-same}. Note that the parameters are assumed to be equal up to first order on $\partial\Omega$, and by our choice of $B$ they are also equal on $\partial B$, so that all boundary integrals on $\partial \Omega_b$ vanish when integrating by parts. 
Thus, the first two terms in \eqref{eq:opp-intformY} become
\begin{equation}
\int_{\Omega_b} U Z^1 \cdot \overline{Y^2} +U Z^1 \cdot \overline{\ddot{Y}^2} dx,  \label{eq:opp-intformU}
\end{equation}
with $U$ as in \eqref{eq:U}. 
We start by computing the limit of the first term. Recalling the asymptotics for the CGO solutions in Proposition \ref{cgo} and the discussion following it, writing
\[ Z^1=e^{i\zeta_1\cdot x}(Z^1_0+Z^1_{-1}+Z^1_r),~~Y^2=e^{i\zeta_2\cdot x}(Y^2_1+Y^2_0+Y^2_r), \]
we have 
\begin{gather}
\| Z^1_{-1}\|_{L^2_\delta}=O(\tau^{-1}),~~\|Z^1_r\|_{L^2_\delta} = O(\tau^{-(1+\delta')}), \label{eq:opp-asyZ}\\
 \|Y^2_1\|_{L^2(O)}=O(\tau),~~\| Y^2_0\|_{L^2(O)}=O(1),~~\| Y^2_r\|_{L^2(O)}=O(\tau^{-\delta'}) \label{eq:opp-asyY}
 \end{gather}
for sufficiently small $\delta'>0$, so that (recall also that $i\zeta_1+\overline{i\zeta}_2=i\xi$)
\begin{align}
\lim_{\tau\rightarrow\infty}\int_{\Omega_b} \td Y^1\cdot \overline{Y^2}dx &= \lim_{\tau\rightarrow\infty}\int_{\Omega_b} Ue^{i\zeta_1\cdot x}(Z^1_0+Z^1_{-1}+Z^1_r) \cdot \overline{e^{i\zeta_2\cdot x}(Y^2_1+Y^2_0+Y^2_r)} dx \notag \\
&= \lim_{\tau\rightarrow\infty} \int_{\Omega_b} e^{i\xi\cdot x} U Z^1_0\cdot \overline{Y^2_1} dx
+ \lim_{\tau\rightarrow\infty} \int_{\Omega_b} e^{i\xi\cdot x} U Z^1_0\cdot \overline{Y^2_0} dx \notag \\
&~~~~~~~+ \lim_{\tau\rightarrow\infty} \int_{\Omega_b} e^{i\xi\cdot x} U Z^1_{-1}\cdot \overline{Y^2_1} dx. \label{eq:opp-limitsleft}
\end{align}
These are the same terms as were obtained in the proof of Theorem \ref{thm-same}, and the limits are summarized in Lemma \ref{same-limlem}.

We proceed by showing that the remaining terms in \eqref{eq:opp-intformY} vanish as $\tau$ becomes large. In order to do so, we first express the vectors $\zeta_1$ and $\zeta_2$ in a different basis, which will allow us to efficiently compute the exponentials resulting from multiplying the reflected solutions. As was done in \cite{COS2009}, we introduce the orthonormal basis $\{f_1,f_2,f_3\}$ of $\RR^3$ with
\begin{equation}
f_2=\frac{1}{|\xi'|}(\xi_1,\xi_2,0),~~~f_3=(0,0,1)=e_3,~~~f_1=f_2\wedge f_3.
\label{eq:opp-onbf}
\end{equation}
With respect to this basis, we have
\[ \eta_1=(1,0,0)_f=f_1,~~~\eta_2=|\xi|^{-1}(0,-\xi_3,|\xi'|)_f,~~~\xi=(0,|\xi'|,\xi_3)_f. \]
Then, the vectors $\zeta_1$ and $\zeta_2$ are
\begin{gather}
\zeta_1=\left( \sqrt{\tau^2-\frac{|\xi|^2}{4}}, \frac{|\xi'|}{2} + i \sqrt{\tau^2-k^2},\frac{\xi_3}{2}-i\sqrt{\tau^2-k^2}\right)_f,\\
\zeta_2=\left( \sqrt{\tau^2-\frac{|\xi|^2}{4}}, -\frac{|\xi'|}{2} - i \sqrt{\tau^2-k^2},-\frac{\xi_3}{2}+i\sqrt{\tau^2-k^2}\right)_f.
\end{gather}
Clearly, for two points $x=(x_1,x_2,x_3)=(x_1^f,x_2^f,x_3^f)_f$ and $y=(y_1,y_2,y_3)=(y_1^f,y_2^f,y_3^f)_f$ in $ \RR^3$ we have $x_3=x_3^f$ and $y_3=y_3^f$, and
\[ x\cdot y = \sum_{j=1}^3x_jy_j =\sum_{j=1}^3x_j^fy_j^f. \]
In terms of this basis, the exponents in the last three terms in \eqref{eq:opp-intformY} are easily computed as
\begin{gather}
\varphi_2(x):= i\zeta_1\cdot x + \overline{i\zeta_2}\cdot \ddot{x} = i |\xi'| x_2^f + i\xi_3 L -2\sqrt{\tau^2-k^2}\frac{|\xi'|}{|\xi|} (L-x_3^f),\label{eq:opp-exp2}\\
\varphi_3(x):=i\zeta_1\cdot \dot{x} + \overline{i\zeta_2}\cdot x = i |\xi'| x_2^f + i\xi_3 L -2\sqrt{\tau^2-k^2}\frac{|\xi'|}{|\xi|} x_3^f, \label{eq:opp-exp3} \\
\varphi_4(x):=i\zeta_1\cdot\dot{x}+\overline{i\zeta_2}\cdot\ddot{x} = i\dot{\xi}\cdot x + i\xi_3 L - 2 \sqrt{\tau^2-k^2}\frac{|\xi'|}{|\xi|} L,\label{eq:opp-exp4}
\end{gather} 
where $\dot{\xi}=(\xi_1,\xi_2,-\xi_3)$. Note that all the exponentials $e^{\varphi_j(x)}$ decay to zero as $\tau\rightarrow \infty$.

\subsubsection*{Limit of the second term in \eqref{eq:opp-intformY}}
We already rewrote this term as in \eqref{eq:opp-intformU}. { Note that this is the analog of the integral \eqref{eq:same-lim4} in the previous section. We proceed to manipulate the term in a fashion similar to \cite[Lemma 4.4]{COS2009}, but our final argument in showing that this integral vanishes will be different due to the different choice of phase vectors. We will employ the decay of the resulting exponential, as was done in \cite{KLU2012}.} 

Recall that the resulting exponent after substituting the solutions is \eqref{eq:opp-exp2}, and we can write the exponential as (using the fact that $x_3=x_3^f$)
\[ e^{\varphi_2(x)} 
= \frac{1}{2\sqrt{\tau^2-k^2}\frac{|\xi'|}{|\xi|}} \partial_{x_3} e^{\varphi_2(x)}, \]
so that we get 
\begin{align} 
\int_{\Omega_b} U Z^1 \cdot \overline{\ddot{Y}^2} dx &=  \frac{1}{2\sqrt{\tau^2-k^2}\frac{|\xi'|}{|\xi|}}\int_{\Omega_b} \partial_{x_3} e^{\varphi_2(x)} U (Z_{-1}^1+Z_0^1+Z_r^1)\cdot \overline{(\ddot{Y}^2_1+\ddot{Y}^2_0+\ddot{Y}^2_r)} dx \notag \\
&= - \frac{1}{2\sqrt{\tau^2-k^2}\frac{|\xi'|}{|\xi|}}\int_{\Omega_b}  e^{\varphi_2(x)} \partial_{x_3} \big[U (Z_{-1}^1+Z_0^1+Z_r^1)\cdot \overline{(\ddot{Y}^2_1+\ddot{Y}^2_0+\ddot{Y}^2_r)} \big] dx \label{eq:opp-lim2}\\
&\qquad \qquad +  \frac{1}{2\sqrt{\tau^2-k^2}\frac{|\xi'|}{|\xi|}}\int_{\partial \Omega_b} \nu_3 e^{\varphi_2(x)} U (Z_{-1}^1+Z_0^1+Z_r^1)\cdot \overline{(\ddot{Y}^2_1+\ddot{Y}^2_0+\ddot{Y}^2_r)} dx.\notag
\end{align}
Now the boundary term vanishes by our assumptions on the parameters. The partial derivative yields the terms 
\begin{align*}
\partial_{x_3} \big[U (Z_{-1}^1+Z_0^1+Z_r^1)\cdot \overline{(\ddot{Y}^2_1+\ddot{Y}^2_0+\ddot{Y}^2_r)} \big] =&~ \partial_{x_3} U (Z_{-1}^1+Z_0^1+Z_r^1)\cdot \overline{(\ddot{Y}^2_1+\ddot{Y}^2_0+\ddot{Y}^2_r)} \\
& ~~+ U \partial_{x_3}(Z_{-1}^1+Z_r^1)\cdot \overline{(\ddot{Y}^2_1+\ddot{Y}^2_0+\ddot{Y}^2_r)} \\
& ~~+ U  (Z_{-1}^1+Z_0^1+Z_r^1)\cdot\partial_{x_3}  \overline{(\ddot{Y}^2_0+\ddot{Y}^2_r)},
\end{align*}
since $Z^1_0$ and $\ddot{Y}^2_1$ are constant with respect to $x$. From the asymptotic behavior \eqref{eq:opp-asyZ}-\eqref{eq:opp-asyY} it follows that
\begin{align}
\partial_{x_3} U (Z_{-1}^1+Z_0^1+Z_r^1)\cdot \overline{(\ddot{Y}^2_1+\ddot{Y}^2_0+\ddot{Y}^2_r)} &=\partial_{x_3} U Z_0^1\cdot \overline{\ddot{Y}^2_1} + O(1),\label{eq:opp-partials1}\\
U \partial_{x_3}(Z_{-1}^1+Z_r^1)\cdot \overline{(\ddot{Y}^2_1+\ddot{Y}^2_0+\ddot{Y}^2_r)} &= U \partial_{x_3}Z_r^1\cdot \overline{\ddot{Y}^2_1} +O(1),\label{eq:opp-partials2}\\
U  (Z_{-1}^1+Z_0^1+Z_r^1)\cdot\partial_{x_3}  \overline{(\ddot{Y}^2_0+\ddot{Y}^2_r)} &= UZ_0^1\cdot\partial_{x_3}  \overline{\ddot{Y}^2_r} +O(1)\label{eq:opp-partials3}
\end{align}
in $L^2(\Omega_b)$. Furthermore, the terms in \eqref{eq:opp-partials2} and \eqref{eq:opp-partials3} satisfy
\[ \|U \partial_{x_3}Z_r^1\cdot \overline{\ddot{Y}^2_1}\|_{L^2(\Omega_b)} = O(\tau^{1-\delta'}),~~~  \| UZ_0^1\cdot\partial_{x_3}  \overline{\ddot{Y}^2_r} \|_{L^2(\Omega_b)} = O(\tau^{1-\delta'}), \]
so that as $\tau\rightarrow\infty$ in \eqref{eq:opp-lim2}, we are left only with the limit of \eqref{eq:opp-partials1},
\[ \lim_{\tau\rightarrow \infty} \int_{\Omega_b} U Z^1 \cdot \overline{\ddot{Y}^2} dx = - \lim_{\tau\rightarrow\infty} \frac{1}{2\sqrt{\tau^2-k^2}\frac{|\xi'|}{|\xi|}}\int_{\partial \Omega_b} e^{\varphi_2(x)} \partial_3(U) Z^1_0 \cdot \overline{\ddot{Y}^2_1} dx.\]
Using the definitions of $Z^1_0$ and $Y^2_1$, 
\[ Z^1_0 = \frac{1}{\tau}\big( \zeta_1\cdot a_1,0, \zeta_1\cdot b_1,0\big) + O(\tau^{-1}),~~~Y^2_1= -\frac{1}{\tau}\big(0, (\zeta_2\cdot b_2)\zeta_2,0,(\zeta_2\cdot a_2)\zeta_2\big) + O(1), \]
and that of $U$, we find 
\begin{multline*} 
- \lim_{\tau\rightarrow\infty} \frac{1}{2\sqrt{\tau^2-k^2}\frac{|\xi'|}{|\xi|}}\int_{\partial \Omega_b} e^{\varphi_2(x)} \partial_3(U) Z^1_0 \cdot \overline{\ddot{Y}^2_1} dx  \\ 
= 2i\omega \int_{\Omega_b} \frac{1}{2 \tau^2 \sqrt{\tau^2-k^2}\frac{|\xi'|}{|\xi|}} e^{\varphi_2(x)} \Big[ (\zeta_1\cdot b_1)(\overline{\zeta_2\cdot b_2}) (\overline{\zeta_2}\cdot \nabla) \partial_{x_3} \Big(\frac{\gamma_1}{\gamma_2}\Big)^{1/2} + (\zeta_1\cdot a_1)(\overline{\zeta_2\cdot a_2}) (\overline{\zeta_2}\cdot \nabla) \partial_{x_3} \Big(\frac{\mu_1}{\mu_2}\Big)^{1/2}\Big] dx.
\end{multline*}
Noting that for $\tau>k$ we have $|\zeta_j/\tau|^2 =2-k^2/\tau^2 \leq 2$, the integrand is bounded in absolute value by
\[ C(|\xi|,|a_1|,|a_2|,|b_1|,|b_2|)\bigg( \Big|\nabla \partial_{x_3} \Big(\frac{\gamma_1}{\gamma_2}\Big)^{1/2}\Big| + \Big|\nabla \partial_{x_3} \Big(\frac{\mu_1}{\mu_2}\Big)^{1/2}\Big| \bigg),  \]
which is an integrable function on $\Omega_b$, since the material parameters are sufficiently smooth bounded functions. Furthermore, the integrand decays to zero pointwise as $\tau \rightarrow \infty$, so by the Dominated Convergence Theorem, the integral vanishes in the limit. This shows that the second term in \eqref{eq:opp-intformY} vanishes as $\tau\rightarrow \infty$.

\subsubsection*{Limit of the third term in \eqref{eq:opp-intformY}}
The third term can be rewritten similarly to the first two, by substituting $Y^2=(P(i\nabla)+\check{W}_2^T)Z^2$. Then, using arguments like those for rewriting the first and second terms, we arrive at 
\[ \int_{\Omega_b} \td \dot{Y}^1\cdot \overline{Y^2}dx = \int_{\Omega_b} \dot{Y}^1\cdot \overline{\tilde{U} Z^2}dx, \]
where
\[ \tilde{U}=\overline{\kappa_2} \tdo + \overline{\kappa_1} \tdgo - i  \left(\begin{array}{cccc}
0 &0&0&\nabla \hat{\mu}\cdot\\
0 & 0 &\nabla \hat{\mu} &\nabla \hat{\mu}\wedge \\
0 &\nabla \overline{\hat{\gamma}}\cdot & 0 & 0 \\
\nabla \overline{\hat{\gamma}}&-\nabla \overline{\hat{\gamma}}\wedge & 0 & 0 \end{array} \right) - i  \left(\begin{array}{cccc}
0 &0&0&\nabla \tilde{\mu}\cdot\\
0 & 0 &\nabla \tilde{\mu} &-\nabla \tilde{\mu}\wedge \\
0 &\nabla \overline{\tilde{\gamma}}\cdot & 0 & 0 \\
\nabla \overline{\tilde{\gamma}}&\nabla \overline{\tilde{\gamma}}\wedge & 0 & 0 \end{array} \right).  \]
The exponent resulting from plugging in the CGO solutions is \eqref{eq:opp-exp3}, and we find that 
\[ e^{\varphi_2(x)} = -\frac{1}{2\sqrt{\tau^2-k^2}\frac{|\xi'|}{|\xi|}} \partial_{x_3} e^{\varphi_2(x)}, \]
so that we can treat this integral like the second term in \eqref{eq:opp-intformY} to find that this integral also vanishes as $\tau \rightarrow \infty$. 

\subsubsection*{Limit of the last term in \eqref{eq:opp-intformY}} 
We finally consider the last term in \eqref{eq:opp-intformY}. With \eqref{eq:opp-exp4}, we obtain 
\begin{gather*}
\int_{\Omega_b} \td \dot{Y}^1\cdot \overline{\ddot{Y}^2}dx = \int_{\Omega_b}\td e^{i\dot{\xi}\cdot x + i\xi_3 L - 2 \sqrt{\tau^2-k^2}\frac{|\xi'|}{|\xi|} L}(\dot{Y}_1^1+\dot{Y}_0^1+\dot{Y}_r^1)\cdot \overline{(\ddot{Y}^2_1+\ddot{Y}^2_0+\ddot{Y}^2_r)}dx\\
\qquad\qquad=e^{ i\xi_3 L - 2 \sqrt{\tau^2-k^2}\frac{|\xi'|}{|\xi|} L}
\int_{\Omega_b} \td e^{i\dot{\xi}\cdot x}\big[ 
(\dot{Y}_1^1+\dot{Y}_0^1+\dot{Y}_r^1)\cdot \overline{(\ddot{Y}^2_1+\ddot{Y}^2_0+\ddot{Y}^2_r)}
\big] dx. 
\end{gather*}
By the asymptotics for the CGO solutions, the integral is $O(\tau)$ as $\tau\rightarrow\infty$, so that in the limit, this term vanishes.

Summarizing, the only terms left when taking the limit in \eqref{eq:opp-intformY} are those in \eqref{eq:opp-limitsleft}. These limits are stated in Lemma \ref{same-limlem}, and 
by choosing proper values of $a_1,b_1$ and $a_2,b_2$ and manipulating the resulting terms conveniently, we arrive at a set of partial differential equations.

\newdefi{Proposition} \label{pdes-opp} {\it Let $\xi\in\RR^3$. If we let $a_1=\overline{a_2}=\tilde{\zeta}=\eta_1-i\eta_2$, and $b_1=\overline{b_2}=\overline{\tilde{\zeta}}=\eta_1+i\eta_2$, then 
\[
\int_O\exi \bigg( \frac{1}{2} \nabla \cdot (\beta_2-\beta_1) +\frac{1}{4}\big( \beta_2\cdot\beta_2-\beta_1\cdot \beta_1\big)+ \kappa_1^2-\kappa_2^2\bigg) dx=0.
\]
If we let $a_1=\overline{a_2}=\overline{\tilde{\zeta}}$ and $b_1=\overline{b_2}=\tilde{\zeta}$, then 
\[
\int_O\exi \bigg( \frac{1}{2} \nabla \cdot (\alpha_2-\alpha_1) +\frac{1}{4}\big( \alpha_2\cdot \alpha_2-\alpha_1\cdot \alpha_1\big)+ \kappa_1^2-\kappa_2^2\bigg) dx=0.
\]
Setting $u=(\gamma_1/\gamma_2)^{1/2}$ and $v=(\mu_1/\mu_2)^{1/2}$, we can rewrite the resulting equations as 
\begin{align*}
- \nabla\cdot (\mu_2\nabla v)+ \omega^2\mu_2^2\gamma_2(u^2v^2-1)v&=0,\\
-\nabla \cdot (\gamma_2\nabla u) + \omega^2\mu_2\gamma_2^2(u^2v^2-1)u&=0,
\end{align*}
and we have the following boundary conditions for $u$ and $v$ by our assumptions on the parameters:
\[ u = v = 1,~~ \frac{\partial u}{\partial \nu} = \frac{\partial v}{\partial \nu} =0~~\on \partial \Omega_b. \]
}
The proof is analogous to that of Proposition \ref{pdes}. The uniqueness proof is now completed by referring to Proposition \ref{uniquecont} with $O=\Omega_b$, which allows us to conclude that $u=v\equiv 1$ in $\Omega_b$.

\begin{appendices}
\renewcommand{\thesection}{A}
\renewcommand{\thesubsubsection}{A.\arabic{subsubsection}}
\section{Well-posedness of the direct problem}

In this section we show the unique solvability of \eqref{eq:opp-ME1}-\eqref{eq:opp-ME3} under the stated regularity assumptions. 
We first reduce the problem to one with zero boundary condition by using the trace theorem: let $E_o \in H^1(\mathbb{R}^3)^3$ be a compactly supported function such that $\nu\wedge E_o=f$ on $\Gamma_1$ and $\nu\wedge E_o=0$ on $\Gamma_2$, and look for a pair of solutions $(\tilde{E},\tilde{H})=(E+E_o, H)$, where $E$ and $H$ solve 
\begin{subequations}
 \label{eq:MEzbc}
 \begin{align}
\nabla\wedge E(x) - i\omega\mu(x) H(x) &= -\nabla\wedge E_o(x) =: F_1(x) \label{eq:ME1zbc} \\
\nabla\wedge H(x) + i\omega\gamma(x) E(x)& = -i\omega\gamma(x)E_o(x)  =: F_2(x) ~~~~\mathrm{in}~\Omega,~~~~~ \nu \wedge E\big|_{\partial\Omega} =0, \label{eq:ME2zbc}
\end{align}
\end{subequations} 
and satisfy a suitable radiation condition as $|(x_1,x_2)|\rightarrow \infty$ for $E$ and $H$. Note that $F_1$ and $F_2$ are compactly supported and we have $F_1 \in L^2(\mathbb{R}^3)^3$ and $F_2 \in H^1(\mathbb{R}^3)^3$. 

In order to solve the system \eqref{eq:MEzbc}, we use the Lax Phillips method and split it into two separate problems, one with constant coefficients on the whole slab and one with nonconstant coefficients, but on a suitable bounded domain $\hat{\Omega}$: we fix $R>0$ such that $\mu$ and $\varepsilon$ are constant outside the ball $B(0,R)$, then choose $R''>R$, and let $\hat{\Omega}$ be a bounded convex domain with $C^{2,1}$ boundary in $\Omega$ containing $\Omega\cap B(0,R'')$, such that $\omega$ is not an eigenvalue for the Maxwell system with constant coefficients $\mu_o,\varepsilon_o$ on $\hat{\Omega}$. 

Now we fix $R'$ with $R<R'<R''$ and pick a function $\varphi \in C_c^\infty(\mathbb{R}^3)$ such that $\varphi(x) = 1$ for all $x\in B(0,R')$ and $\varphi(x)=0$ for all $x\in B(0,R'')^c$. We want to obtain a pair of solutions to \eqref{eq:MEzbc} of the form
\begin{equation}
E= E_1 - \varphi(E_1-E_2),~~~~~H= H_1-\varphi(H_1-H_2), \label{eq:ansatz}
\end{equation} 
where $(E_1,H_1)$ solves Maxwell's equations with constant coefficients in $\Omega$,
\begin{subequations}
 \label{eq:MEconst}
 \begin{align}
\nabla\wedge E_1 - i \omega\mu_o H_1 &= \tilde{F}_1, \label{eq:ME1const} \\
\nabla \wedge H_1 + i\omega\varepsilon_o E_1 &= \tilde{F}_2~~~\mathrm{in}~\Omega,  \label{eq:ME2const}\\
 ~~~~~ \nu\wedge  E_1&=0~~~\mathrm{on}~ \partial\Omega, \label{eq:ME3const}\\ 
\mathrm{radiation~condition~for~}E_1,~H_1~&\mathrm{as}~|(x_1,x_2)|\rightarrow \infty \label{eq:ME4const}
\end{align}
\end{subequations} 
and $(E_2,H_2)$ solves the problem with non-constant coefficients in $\hat{\Omega}$,
\begin{subequations}
 \label{eq:MEnonconst}
 \begin{align}
\nabla\wedge E_2 - i\omega\mu H_2 = \tfo,~~&\label{eq:ME1nonconst} \\
\nabla\wedge H_2 + i\omega\gamma E_2 = \tft~~~&\mathrm{in}~\hat{\Omega}, \label{eq:ME2nonconst}\\
\nu \wedge  E_2 = \nu \wedge E_1 ~~~&\mathrm{on}~ \partial\hat{\Omega}, 
\label{eq:ME3nonconst}
\end{align}
\end{subequations} 
and $\tfo \in L^2(\mathbb{R}^3)^3$ and $\tft \in H^1(\mathbb{R}^3)^3$ 
are compactly supported functions that will be determined in the following. Plugging the ansatz \eqref{eq:ansatz} into \eqref{eq:MEzbc}, the equations are satisfied inside $\Omega \cap B(0,R)$ if $\tfo=F_1$ and $\tft=F_2$ in this region, since $\varphi\equiv 1$ there and $(E_2,H_2)$ solve \eqref{eq:MEnonconst}. In the region $\hat{\Omega}\backslash B(0,R)$ the parameters are constant and we obtain
\begin{eqnarray*}
\nabla \wedge E - i\omega\mu_o H &=&  
\tfo - \nabla \varphi \wedge (E_1-E_2),\\
\nabla \wedge H + i \omega \varepsilon_o E &=&  
\tft - \nabla \varphi \wedge (H_1-H_2),
\end{eqnarray*}
so in order for $(E,H)$ to satisfy \eqref{eq:MEzbc}, $\tilde{F}=(\tfo,\tft)$ needs to satisfy
\begin{equation}
\big( I + K \big) \tilde{F} = \left( \begin{array}{c} F_1 \\ F_2 \end{array} \right), \label{eq:LPsyst}
\end{equation}
where the operator $K$ is defined by
\begin{equation}
K \tilde{F}  
= \left( \begin{array}{c} \nabla \varphi \wedge (E_2-E_1) \\ \nabla \varphi \wedge (H_2 - H_1) \end{array} \right), \label{eq:K}
\end{equation}
with $(E_1,H_1)$ being the solution to \eqref{eq:MEconst} and $(E_2,H_2)$ the solution to \eqref{eq:MEnonconst}. 
We need to show that \eqref{eq:LPsyst} is uniquely solvable. We postpone this and first concern ourselves with the solvability of each of the systems \eqref{eq:MEconst} and \eqref{eq:MEnonconst}; we begin with the constant coefficient system \eqref{eq:MEconst}.

\subsubsection{Maxwell's equations with constant coefficients in the slab}
In the case of constant coefficients we can transform Maxwell's equations to obtain a vector Helmholtz equation:  Taking the divergence of the second equation, we obtain the identiy 
\begin{equation}
\nabla \cdot E_1 = \frac{1}{i\omega\varepsilon_o} \nabla \cdot \tft. \label{eq:divE}
\end{equation}
Taking the curl of the first equation and then using the second equation to substitute $\nabla \wedge H_1$, as well as the identity $\nabla\wedge \nabla \wedge = -\Delta + \nabla \nabla \cdot$ and \eqref{eq:divE}, we arrive at 
\begin{equation}
(-\Delta - k^2 ) E_1 = i \omega \mu_o \tft + \nabla \wedge \tfo - \frac{1}{i\omega\varepsilon_o} \nabla \nabla \cdot \tft =: G,\label{eq:vectorHE}
\end{equation}
with $k= \omega\sqrt{ \mu_o \varepsilon_o}$. The right-hand side $G\in H^{-1}(\mathbb{R}^3)$ is compactly supported. Since on $\partial \Omega$, $\nu = \pm e_3$, the boundary condition gives that $E_{1,1}=E_{1,2}=0$ on $\partial \Omega$. 
So we can expand $E_{1,1}$ and $E_{1,2}$ as sine series in the variable $x_3$,
\[ E_{1,j}(x)=\sum_{m=1}^\infty E_m^j(x_1,x_2) \sin\big( \frac{m \pi x_3}{L}\big). \]
We now address the necessary condition to guarantee uniqueness of $E_{1,1}$ and $E_{1,2}$. This is a radiation condition for the coefficients of this expansion. At this point we need the requirement that $k$ is such that $k \neq m\pi/L$ for all $m\in \mathbb{N}$. 
Then there is a unique solution to the Helmholtz equation in the slab whose coefficients in the above sine series expansion satisfy the following conditions: 
\begin{enumerate}[(i)]
\setlength\itemsep{0em}
\item for $m$ such that $k^2-m^2\pi^2/L^2<0$, $E_m^j \in H^1(\mathbb{R}^2)$;
\item for $m$ such that $k^2-m^2\pi^2/L^2>0$, set $k_m=\sqrt{k^2-m^2\pi^2/L^2}$, then $E^j_m$ satisfies, with $x'=(x_1,x_2)$,
\begin{equation}
E_m^j(x') = O(r^{-1/2}),~~~\big( \frac{\partial}{\partial r} - i k_m\big) E_m^j = o(r^{-1/2}),~~~r=|x'| \rightarrow \infty. \tag{RC} \label{eq:rc}
\end{equation}
\end{enumerate}
The latter is a Sommerfeld radiation condition for $E_m^j$; this so-called partial radiation condition was introduced in \cite{S1950}, and uniqueness of solutions to the Helmholtz equation satisfying this condition was proved in \cite{S1950,RW1985}. We introduce the following notion.\vspace{2mm}

\newdefi{Definition} \label{def-admissible} A solution to \eqref{eq:opp-ME1}-\eqref{eq:opp-ME2} is called \emph{admissible}, if its coefficients in the series expansion satisfy the partial radiation condition \eqref{eq:rc}.\vspace{2mm}

In order to obtain the unique admissible solutions to Maxwell's equations, we use the fundamental solution for the Helmholtz equation in the slab with homogeneous Dirichlet boundary condition,
\begin{equation}
 \Phi(x,y) = \sum_{m=1}^\infty \frac{-1}{2L} \sin\Big(\frac{m\pi x_3}{L}\Big) \sin\Big(\frac{m\pi y_3}{L}\Big) H^1_0(k_m |x'-y'|), \label{eq:FSdirichlet}
\end{equation}
where $k_m = \sqrt{k^2-m^2\pi^2/L^2}$, and $H^1_0$ is the Hankel function of first kind. 
Thus, we obtain the solutions in $H^1(\Omega)$
\begin{equation}
E_{1,1}(x) = \int_\Omega \Phi(x,y) G_1(y) dy, ~~~~E_{1,2}(x) = \int_\Omega \Phi(x,y) G_2(y) dy. \label{eq:sol12}
\end{equation}
This approach is not applicable to $E_{1,3}$, however, since we do not know the boundary value of $E_{1,3}$. To get around this, we use \eqref{eq:divE}, which gives
\[ \partial_{x_3} E_{1,3} = \frac{\partial E_{1,3}}{\partial \nu} = \frac{1}{i\omega\varepsilon_o}\nabla\cdot \tft~~~\mathrm{on}~\partial \Omega. \]
This provides a Neumann boundary condition for $E_{1,3}$, and in order to solve the Helmholtz equation for $E_{1,3}$, we employ the fundamental solution for Neumann boundary data,
\begin{equation}
 \Psi(x,y) = \sum_{m=1}^\infty \frac{-1}{2L} \cos\Big(\frac{m\pi x_3}{L}\Big) \cos\Big(\frac{m\pi y_3}{L}\Big) H^1_0(k_m |x'-y'|).\label{eq:FSneumann}
\end{equation}
We first use the trace theorem again to transform the problem into one with zero boundary condition: Let $\check{E} \in H^1(\mathbb{R}^3)$ with compact support be such that $\partial \check{E}/\partial \nu =1/(i\omega\varepsilon_o)\nabla\cdot \tft$ on $\partial \Omega$, and look for $E_{1,3} = \check{E} + \hat{E}$, where $\hat{E}$ now satisfies
\[ (-\Delta -k^2)\hat{E} = G_3 + (\Delta + k^2)\check{E}~~\mathrm{in}~\Omega,~~~~\frac{\partial \hat{E}}{\partial \nu} = 0 ~~\mathrm{on}~\partial \Omega, \]
as well as a partial radiation condition \eqref{eq:rc} for the coefficients of its cosine expansion. 
Using the fundamental solution $\Psi$ for this problem, we get
\[ \hat{E}(x) = \int_{\Omega} \Psi(x,y) \big( G_3 + (\Delta + k^2)\check{E} \big) dy, \]
and thus
\begin{equation}
E_{1,3}(x)= \check{E}(x) + \int_{\Omega} \Psi(x,y) \big( G_3 + (\Delta + k^2)\check{E} \big) dy.
\label{eq:sol3}
\end{equation}
In order to verify that $E_1=(E_{1,1},E_{1,2},E_{1,3})$ found in this process is indeed a suitable solution for Maxwell's equations, we compute $\nabla \cdot E_1$ to confirm that \eqref{eq:divE} is satisfied. Using the properties of the fundamental solutions $\Phi$ and $\Psi$, this is readily verified. Thus, $E_1\in H^1_{loc}(\Omega)^3$ and $H_1 := \frac{1}{i\omega\mu_o}(\nabla \wedge E_1 -\tfo)$ are solutions to \eqref{eq:ME1const}-\eqref{eq:ME4const}. Note that $H_1\in L^2_{loc}(\Omega)^3$, and by equation \eqref{eq:ME2const}, $\nabla \wedge H_1 \in H^1_{loc}(\Omega)^3$. 

\subsubsection{Maxwell's equations with nonconstant coefficients on $\hat{\Omega}$ and regularity of solutions}
Now we consider the system \eqref{eq:MEnonconst} on the bounded domain $\hat{\Omega}$. Note that since $E_1|_{ \hat{\Omega}} \in H^1(\hat{\Omega})^3$, its tangential trace belongs to $TH^{1/2}(\partial\hat{\Omega})$. 
Under these conditions the system \eqref{eq:MEnonconst} has a unique pair of solutions $(E_2,H_2) \in \Hcurl(\hat{\Omega})^2$, where $\Hcurl(\hat{\Omega}) = \{ u \in L^2(\hat{\Omega})^3: \nabla \wedge u \in L^2(\hat{\Omega})^3\}$. The solutions in fact exhibit higher regularity, as we can verify by using results from \cite[Sec. I.3]{GR-FEM}. 
{ Using \eqref{eq:ME2nonconst}, and the identity \eqref{eq:LPsyst} for $\tft$, we find
\[ \nabla \cdot E_2 = \frac{1}{i\omega \gamma} \big[ \nabla \cdot F_2 - \nabla \varphi \cdot \big(\nabla \wedge (H_1-H_2)\big) - i\omega \nabla \gamma \cdot E_2\big] \in L^2(\hat{\Omega}). \]
Since on $\partial \hat{\Omega}$ we have $\nu\wedge E_2=\nu\wedge E_1 \in TH^{1/2}(\partial\hat{\Omega})$, we can conclude that $E_2 \in H^1(\hat{\Omega})$. Similarly, using \eqref{eq:ME1nonconst},
\[ \nabla \cdot H_2 = \frac{1}{i\omega\mu} \big[\nabla \varphi \cdot \big( \nabla \wedge (E_1-E_2) \big) - i\omega \nabla \mu \cdot H_2 \big] \in L^2(\hat{\Omega}).\]
Finally for $H_1$, taking the divergence in \eqref{eq:ME1const} using the definition of $\tfo$ as well as that of $F_1$, we compute
\[ \nabla \cdot H_1 = -\frac{1}{i \omega\mu_o} \nabla \cdot \tfo = -\frac{1}{i \omega\mu_o} \nabla \cdot ( F_1 - \nabla \varphi \wedge (E_2-E_1)) = \frac{1}{i \omega\mu_o} \nabla \varphi \cdot \big(\nabla \wedge (E_2-E_1) \big) \in L^2(\Omega),  \]
and on $\partial \Omega$ we have $i\omega\mu_o \nu \cdot H_1 = \nabla_{\partial \Omega} \cdot f \in H^{1/2}(\partial \Omega)$, so that summarizing, we have}
\[ H_1 \in L^2_{loc}(\Omega)^3,~~\nabla\wedge H_1 \in H^1_{loc}(\Omega),~~\nabla\cdot H_1 \in L^2(\Omega),~~\nu\cdot H_1 \in H^{1/2}(\partial \Omega). \]

\subsubsection{Compactness of $K$ and unique solvability of \eqref{eq:LPsyst}}
We proceed to show that the operator $K: L^2(\hat{\Omega})^3\times H^1(\hat{\Omega})^3 \rightarrow L^2(\hat{\Omega})^3\times H^1(\hat{\Omega})^3$ defined in \eqref{eq:K} is indeed compact. We will do so by showing that $KF \in H^1(\hat{\Omega})^3\times H^2(\hat{\Omega})^3$ and using compact embedding. 

We first note that since $E_1|_{\hat{\Omega}}$ and $E_2$ belong to $H^1(\hat{\Omega})^3$ and $\varphi \in C_c^\infty(\mathbb{R}^3)$, $\nabla \varphi \wedge (E_2-E_1) \in H^1(\hat{\Omega})^3$. It remains to show that $\nabla \varphi \wedge (H_2-H_1) \in H^2(\hat{\Omega})^3$. 
Note that since wherever $\nabla \varphi \neq 0$, we have constant parameters $\mu = \mu_o$ and $\varepsilon=\varepsilon_o$, in this region we also have the identities
\begin{align}  
\nabla\cdot (H_2-H_1)&=0, \label{eq:Kpf1} \\
\nabla \wedge (H_2-H_1)&= - i \omega \varepsilon_o (E_2-E_1) \in H^1(\hat{\Omega})^3. \label{eq:Kpf3}
\end{align}

We take a smooth cutoff function $\psi \in C^\infty_c(\mathbb{R}^3)$ that satisfies $\psi(x)=0$ on $B(0,R)$, and $\psi(x) = 1$ on $B(0,R_{\psi,2})\backslash B(0,R_{\psi,1})$ for radii $R< R_{\psi,1}<R'<R''<R_{\psi,2}$, such that $\psi\equiv 1$ wherever $\nabla \varphi\neq 0$. Using \eqref{eq:Kpf1} and \eqref{eq:Kpf3}, we then compute
\begin{gather}
\nabla \cdot \psi (H_2-H_1) =  \nabla \psi \cdot (H_2-H_1) + \psi \nabla \cdot (H_2-H_1) = \nabla \psi \cdot (H_2-H_1) \in L^2(\hat{\Omega}),\label{eq:diffdiv}\\
\nabla\wedge \psi (H_2-H_1) = \nabla \psi \wedge (H_2-H_1) + \psi \nabla \wedge (H_2-H_1)  \in L^2(\hat{\Omega})^3. \label{eq:diffcurl}
\end{gather}
{ We also have $\nu \cdot \psi (H_2-H_1) =0$ on $\partial \hat{\Omega}$, using
\begin{equation*}
i \omega \nu \cdot ( \mu_o H_1) = \nu \cdot ( \nabla \wedge E_1 -\tfo ) = -\nabla_{\partial \hat{\Omega}} \cdot (\nu \wedge E_1) - \nu \cdot \tfo =-\nabla_{\partial \hat{\Omega}} \cdot (\nu \wedge E_2) - \nu \cdot \tfo = i \omega \nu \cdot ( \mu H_2)
\end{equation*}
and the fact that $\mu=\mu_o$ where $\psi \neq 0$. This implies that $\psi (H_2-H_1) \in H^1(\hat{\Omega})^3$ \cite[Cor. I.3.7]{GR-FEM}. Using this fact in \eqref{eq:diffdiv} and \eqref{eq:diffcurl} in turn implies that $\nabla \cdot \psi (H_2-H_1) \in H^1(\hat{\Omega})$ and  $\nabla\wedge \psi (H_2-H_1) \in H^1(\hat{\Omega})^3$, and since $\hat{\Omega}$ is a $C^{2,1}$ domain, we can apply Corollary I.3.7 from \cite{GR-FEM} again to conclude that $\psi(H_2-H_1)\in H^2(\hat{\Omega})$. Thus,
\[ \nabla \varphi \wedge (H_2-H_1) =\nabla \varphi \wedge \psi (H_2-H_1) \in H^2(\hat{\Omega})^3.  \]}

We find that $K$ maps into $H^1(\hat{\Omega})^3\times H^2(\hat{\Omega})^3$, and by compact embedding, $K$ is a compact operator on $L^2(\hat{\Omega})^3\times H^1(\hat{\Omega})^3$.

Now we can use Fredholm theory to show unique solvability for \eqref{eq:LPsyst}. Solvability follows once we have established uniqueness of the solution. So suppose that $F_1=F_2=0$. Then $(E,H)$ solves the homogeneous Maxwell equations with zero tangential boundary condition for $E$; thus $E=H=0$ by Assumption 1.

In $\Omega\cap B(0,R)$, by definition $K\tilde{F}=0$, and hence \eqref{eq:LPsyst} yields $\tilde{F}=(F_1,F_2)=(0,0)$ in this region. Furthermore, \eqref{eq:ansatz} shows $E_2=E=0$ and $H_2=H=0$ in $\Omega\cap B(0,R)$. Thus in particular, $E_2$ and $H_2$ solve Maxwell's equations with constant coefficients $\mu_o$ and $\varepsilon_o$ in $\Omega\cap B(0,R)$. Since outside this region, the coefficients are constant, we conclude that 
\begin{align*} 
\nabla \wedge (E_1-E_2) - i \omega \mu_o (H_1-H_2) &= 0 \\
\nabla \wedge (H_1-H_2) + i\omega \varepsilon_o (E_1-E_2) &= 0~~ \mathrm{in}~\hat{\Omega},~~~~~\nu\wedge (E_1-E_2)=0~~\mathrm{on}~\partial\hat{\Omega}, 
\end{align*} 
and since $\omega$ is not an eigenvalue on $\hat{\Omega}$, this yields $E_1=E_2$ and $H_1=H_2$ in $\hat{\Omega}$, thus $K\tilde{F}=0$ and $\tilde{F}=(F_1,F_2)=(0,0)$ in $\hat{\Omega}$. 

\end{appendices}

\begin{center}
{\sc Acknowledgments}
\end{center}
The author would like to thank her advisor, Ting Zhou, for her guidance and for many helpful discussions. The author's research was supported by the NSF grant DMS-1544138.

\bibliographystyle{alpha} 
\bibliography{inverse.bib}

{\sc Department of Mathematics, Northeastern University. Boston, MA 02115. USA.} 

E-mail address: \href{mailto:pichler.mo@husky.neu.edu}pichler.mo@husky.neu.edu

\end{document}